%% file: Revisedmain.tex
\documentclass[sn-basic,Numbered]{sn-jnl}%

\usepackage{graphicx}
\usepackage{multirow}
\usepackage{amsmath,amssymb,amsfonts}%
\usepackage{amsthm}%
\usepackage{mathrsfs}%
\usepackage[title]{appendix}%
\usepackage[dvipsnames]{xcolor}
\usepackage{textcomp}%
\usepackage{manyfoot}%
\usepackage{booktabs}%
\usepackage{tikz,pgfplots}
\pgfplotsset{compat=1.15}
\usetikzlibrary{shapes,arrows,matrix,fit,shapes,backgrounds}
\usepackage{mathrsfs,wasysym}
\usepackage{siunitx}
\usepackage{algorithm}%
\usepackage{algorithmicx}%
\usepackage{algpseudocode}%
\usepackage{listings}%

\usepackage{subfig}
\hypersetup{
linkcolor=NavyBlue
,citecolor=NavyBlue
,filecolor=Red
,urlcolor=NavyBlue
,menucolor=Red
,runcolor=Red
,linkbordercolor=NavyBlue
,citebordercolor=NavyBlue
,filebordercolor=Red
,urlbordercolor=NavyBlue
,menubordercolor=Red
,runbordercolor=Red
}

\theoremstyle{thmstyleone}%

\theoremstyle{thmstyletwo}%
\newtheorem{remark}{Remark}%

\theoremstyle{thmstylethree}%

\begin{document}

\title[Alya towards Exascale]{Alya towards Exascale: Algorithmic Scalability using PSCToolkit}

\author[1]{\fnm{Herbert} \sur{Owen}}\email{herbert.owen@bsc.es}

\author[1]{\fnm{Oriol} \sur{Lehmkuhl}}\email{oriol.lehmkuhl@bsc.es}

\author*[2]{\fnm{Pasqua} \sur{D'Ambra}}\email{pasqua.dambra@cnr.it}

\author[3,2]{\fnm{Fabio} \sur{Durastante}}\email{fabio.durastante@unipit.it}

\author[4,2]{\fnm{Salvatore} \sur{Filippone}}\email{salvatore.filippone@uniroma2.it}

\affil*[1]{\orgname{Barcelona Supercomputing Centre (BSC)}, \orgaddress{\street{Plaça d'Eusebi Güell}, \city{Barcelona}, \postcode{08034}, \country{Spain}}}

\affil[2]{\orgdiv{Institute for Applied Computing}, \orgname{National Research Council (CNR)}, \orgaddress{\street{Via P. Castellino, 111}, \city{Naples}, \postcode{80131}, \state{NA}, \country{Italy}}}

\affil[3]{\orgdiv{Department of Mathematics}, \orgname{University of Pisa}, \orgaddress{\street{Largo Bruno Pontecorvo, 5}, \city{Pisa}, \postcode{56127}, \state{PI}, \country{Italy}}}

\affil[4]{\orgdiv{Department of Civil and Computer Engineering}, \orgname{University of Rome ``Tor Vergata''}, \orgaddress{\street{Via del Politecnico, 1}, \city{Rome}, \postcode{00133}, \state{RM}, \country{Italy}}}

\abstract{In this paper, we describe an upgrade of the Alya code with up-to-date parallel linear solvers capable of achieving reliability, efficiency and scalability in the computation of the pressure field at each time step of the
numerical procedure for solving a {Large Eddy Simulation} formulation of the incompressible Navier-Stokes equations. We developed a software module in Alya's kernel to interface the libraries included in the current version of \texttt{PSCToolkit}, a framework
for the iterative solution of sparse linear systems on parallel distributed-memory computers by Krylov methods coupled to Algebraic MultiGrid preconditioners. The Toolkit has undergone various extensions within the EoCoE-II project with the primary goal of facing the exascale challenge.
Results on a realistic benchmark for airflow simulations in wind farm applications show that the \texttt{PSCToolkit} solvers
significantly outperform the original versions of the Conjugate Gradient method available in the Alya's kernel in terms of scalability and parallel efficiency and represent a very promising software layer to move the Alya code towards exascale.}



\maketitle

\section{Introduction} \label{Introduction}

Alya is a high-performance computational mechanics code for complex coupled multi-physics engineering problems. In this work, we present the interfacing between Alya and the \texttt{PSCToolkit} to overcome one of Alya's main obstacles in the path towards exascale, namely the lack of state-of-the-art parallel algebraic linear solvers with adequate algorithmic scalability, as already identified in~\cite{VAZQUEZ201615}, where Alya's strengths and weaknesses in facing the exascale challenge have been analyzed by scalability studies up to one hundred thousand cores.

Although Alya can be applied to a wide range of problems, in this work, we shall concentrate {on solving turbulent incompressible flow problems using a Large Eddy Simulation (LES) approach}. Due to the wide range of scales present in turbulent high-Reynolds-number flows, their accurate solution requires computational meshes with a huge number of degrees of freedom (dofs). %
Alya uses a Finite Element (FE) spatial discretization, while its time discretization is based on finite difference methods; when an implicit time discretization is applied, the two main kernels of a simulation are the assembly of stiffness matrices and the solution of the associated linear system at each time step. {In~\cite{VAZQUEZ201615} the authors observed that the FE assembly implemented in Alya showed nearly perfect scalability, as one could a priory expect, while the solution of linear systems by available iterative linear solvers was the main weakness in the path towards exascale.} The problem is related to Alya's lack of solvers with optimal algorithmic scalability, i.e., solvers able to obtain a given accuracy employing an almost constant number of iterations for an increasing number of dofs.

Alya's sparse linear algebra solvers are specifically developed with tight integration with the overall parallelization scheme; they include {Krylov-based solvers, such as Generalized Minimal Residual (GMRES) or Conjugate Gradient (CG), coupled to some deflation approach or a simple diagonal preconditioner. As shown in~\cite{VAZQUEZ201615}, when incompressible flow problems are considered, the solution of a Poisson-type equation
for the pressure field becomes challenging as the size of the problem increases. Indeed, when a uniform mesh multiplication~\cite{HOUZEAUX2013142} is used to have successively finer mesh,  each time the mesh is refined to obtain elements with half the size, the number of iterations for solving the pressure equation is approximately doubled, showing a mesh-size-dependent behavior. To overcome these scalability issues, we interfaced \texttt{PSCToolkit} to Alya to take advantage of the Algebraic MultiGrid (AMG) preconditioners available through the \texttt{AMG4PSBLAS} library;  this effort has been carried out in the context of the European Center of Excellence for Energy (EoCoE) applications.}

The rest of the paper is organized as follows. In Section~\ref{Alya Description}, we describe the general framework of Alya and the type of fluid dynamics problem we wish to test the new solvers on;  in Section~\ref{sec:PSCToolkit}, we give an overall description of \texttt{PSCToolkit}, and then we focus on the AMG preconditioners employed in Section~\ref{sec:amgprecond}. Section~\ref{interface} discusses the  new module written to interface the solver library to the Alya software and the related issues. Section~\ref{bolund} describes the actual test case, while Section~\ref{results} analyzes the numerical scalability results in detail. Finally, Section~\ref{sec:conclusions} summarizes the results obtained and illustrates the new lines of development. 

\section{Alya Description} \label{Alya Description}

Alya is a high-performance computational mechanics code for complex coupled multi-physics engineering problems. It can solve problems in the simulation of turbulent incompressible/compressible flows, non-linear solid mechanics, chemistry, particle transport, heat transfer, and electrical propagation.
{Alya has been designed for massively parallel supercomputers and exploits several parallel programming models/tools. It relies on MPI to support a distributed-memory model;
some kernels support vectorization at the CPU level and GPU accelerators are exploited through OpenACC pragmas or CUDA.}

Multi-physics coupling is achieved following a multi-code strategy that uses MPI to communicate different instances of Alya. Each instance solves a particular physics, enabling asynchronous execution. Coupled problems can be solved by retaining the scalability properties of the individual instances.
{Alya is one of the two Computational Fluid Dynamic (CFD) codes} of the Unified European Applications Benchmark Suite (UEBAS)~\cite{UEABS}. It is also part of the Partnership for Advanced Computing in Europe (PRACE) Accelerator benchmark suite~\cite{PraceBS}.

As  mentioned in Section~\ref{Introduction}, large-scale CFD applications are the main problems targeted by Alya; hence, 
the  basic mathematical models include various formulations
of the Navier-Stokes equations, whose strong form for 
incompressible flows in a suitable domain is the following:
\begin{align}\label{eq:navierstokes}
\partial_t \mathbf{u} - 2 \nu \nabla \cdot \mathbf{\varepsilon}(\mathbf{u}) + \mathbf{u} \cdot \nabla \mathbf{u} + \nabla p &= \mathbf{f}, \\
\nabla \cdot \mathbf{u} &= 0,
\end{align}

\noindent where $\mathbf{u}$ and $p$ are the velocity and pressure field 
respectively, $\mathbf{\varepsilon}\mathbf{\left(u\right)}=\frac{1}
{2}\left(\nabla\mathbf{u}+\nabla^{T}\mathbf{u}\right)$ is the velocity 
strain rate tensor, $\nu$ is the kinematic viscosity, and $\mathbf{f}$ 
denotes the vector of external body forces. The problem is supplied with an 
initial divergence-free velocity field and appropriate boundary conditions.

The flow is turbulent for most real-world flow problems, and some turbulence modeling is needed to make the problem solvable with currently available computational resources. {For all the examples presented in this work, we rely on the functionalities of Alya, which apply the spatially filtered Navier-Stokes equations coupled to the Vreman subgrid-scale model~\cite{vreman} for turbulence closure. In practice, a spatially varying turbulent viscosity supplements the laminar viscosity and the velocity and pressure unknowns correspond to spatially filtered values.}
Finally, since the size of the dynamically important eddies at high Reynolds numbers becomes too small to be grid resolved close to the wall, we employ a wall modeling technique~\cite{https://doi.org/10.1002/fld.4770} to impose the boundary conditions for the LES equations.
For simplicity, the non-linear term has been written in its convective form, which is most commonly encountered in computational practice. 

{Space discretization is based on a Galerkin FE approximation, employing hybrid unstructured meshes, which can include tetrahedra, prisms, hexahedra, and pyramids.} Temporal discretization is performed through an explicit third-order Runge--Kutta scheme, where the Courant--Friedrichs--Lewy number is set to $\text{CFL} = 1.0$ for the cases presented in this work. 
A non-incremental fractional step method is used to stabilize 
the pressure, allowing the use of {finite 
element} pairs that do not satisfy the inf-sup 
condition~\cite{codinafs}, such as the equal order 
interpolation for the velocity and pressure applied in this 
work. 
A detailed description of the above numerical method, together 
with examples for turbulent flows, showing its high accuracy 
and low dissipation, can be found in~\cite{Lehmkuhl_low_dis}.

The fractional step method allows uncoupling the solution of velocity and pressure~\cite{codinafs}. At each Runge--Kutta substep, {an explicit approach computes an intermediate velocity, and }then a linear system coming from a Poisson-type equation is solved for the pressure; 
finally, the incompressible velocity is recovered. In the path towards exascale, the solution of the linear system for the pressure is the most demanding step.  
{To} reduce the computational burden, for most problems, an approximate projection method for Runge--Kutta time-stepping schemes is applied, {which}  allows solving for the pressure only at the final substep~\cite{doi:10.2514/1.J054569}.

It is important to note that most flow problems solved with Alya use a fixed mesh. For such problems, the linear system matrix for the pressure equation remains constant during the whole simulation. Therefore,
the matrix assembly and the setup of a matrix preconditioner are needed only once at the beginning of the numerical procedure. Given that the number of time steps for LES is usually of the order of $10^5$, it is clear that the linear solver computational times and scalability are the most relevant issues to be tackled. 

\section{\texttt{PSCToolkit}: \texttt{PSBLAS} and \texttt{AMG4PSBLAS}}\label{sec:PSCToolkit}

We have interfaced Alya to exploit the solvers and preconditioners developed in the \texttt{PSCToolkit}\footnote{See~\href{https://psctoolkit.github.io/}{psctoolkit.github.io}
on how to obtain and run the code.} software framework for parallel sparse 
computations, proven on current petascale supercomputers and targeting the next-generation exascale machines.
\texttt{PSCToolkit} is composed of two main libraries, named \texttt{PSBLAS} (Parallel
Sparse Basic Linear Algebra Subprograms)~\cite{FC2000,FB2012}, and
\texttt{AMG4PSBLAS} (Algebraic MultiGrid Preconditioners for
PSBLAS)~\cite{dambra2020b}.

Both libraries are written in modern Fortran; \texttt{PSBLAS} implements 
algorithms and functionalities of parallel iterative Krylov subspace linear 
solvers, while  \texttt{AMG4PSBLAS} is the package containing sophisticated 
preconditioners. In particular, \texttt{AMG4PSBLAS} provides 
one-level Additive Schwarz (AS) and Algebraic MultiGrid (AMG) preconditioners. In
the following, we will describe in some detail the AMG preconditioners 
we use within the Alya test cases.

\subsection{AMG preconditioners}\label{sec:amgprecond}

Algebraic MultiGrid methods can be viewed as a particular instance of a general stationary iterative method:
\[
\mathbf{x}^{(k)}=\mathbf{x}^{(k-1)}+B\left(\mathbf{b}-A\mathbf{x}^{(k-1)}\right), \ \ k=1, 2, \ldots\; \text{given}\; \mathbf{x}^{(0)} \in \mathbb{R}^n,
\]
for the solution of a linear system 
\[
A \mathbf{x} = \mathbf{b}, \qquad A \in \mathbb{R}^{n \times n}, \; \mathbf{b} \in \mathbb{R}^n,
\]
where $A$ is symmetric and positive-definite (SPD), and the iteration matrix 
$B$ is defined recursively; see, e.g., \cite{MR2427040} for an exhaustive 
account. 
AMG methods are often employed as preconditioners for Krylov subspace 
solvers; what distinguishes the methods implemented in \texttt{AMG4PSBLAS} 
are the  specific details of the construction procedure for the $B$ matrix. 

We  define $A_0 = A$, and consider the sequence $\{A_l\}_{l=0}^{n_\ell-1}$ of coarse matrices computed by the triple-matrix Galerkin product:
\[
A_{l+1}=P_l^T A_l P_l, \ \ l=0, \ldots, n_\ell-1,
\]
where $\{ P_l \}_{l = 0}^{n_\ell-1}$ is a sequence of prolongation matrices of size $n_l \times n_{l+1}$, with $n_{l+1} < n_l$ and $n_0=n$. To complete the formal construction we need also a sequence $\{ M_l \}_{l = 0}^{n_\ell-2}$ of $A_l$-convergent smoothers for the coarse matrices $A_l$, i.e., matrices $M_l$ for which  $\|I_l -M_l^{-1}A_l\|_{A_l} < 1$ holds true, where $I_l$ is the identity matrix of size $n_l$ and $\| \mathbf{v} \|_{A_l} = \sqrt{\mathbf{v}^T A_l \mathbf{v}}$ is the $A_l$ norm. The preconditioner matrix $B$ for the $V$-cycle with $\nu$ pre- and post-smooth iteration is then given by the multiplicative composition of the following error propagation matrices,
\begin{equation}\label{eq:multiplicative_composition}
    I_l -B_lA_l= (I_l-M_l^{-T}A_l)^\nu(I_l-P_lB_{l+1}P_l^TA_l)(I_l-M_l^{-1}A_l)^\nu \, \forall l < n_{\ell}, 
\end{equation}
with $B_{n_{\ell}} \approx A_{n_\ell}^{-1}$, either as a \emph{direct solution} or as a convergent iterative procedure with a fine enough tolerance.

\begin{figure}[htbp]
    \centering
    \includegraphics[width=0.65\columnwidth]{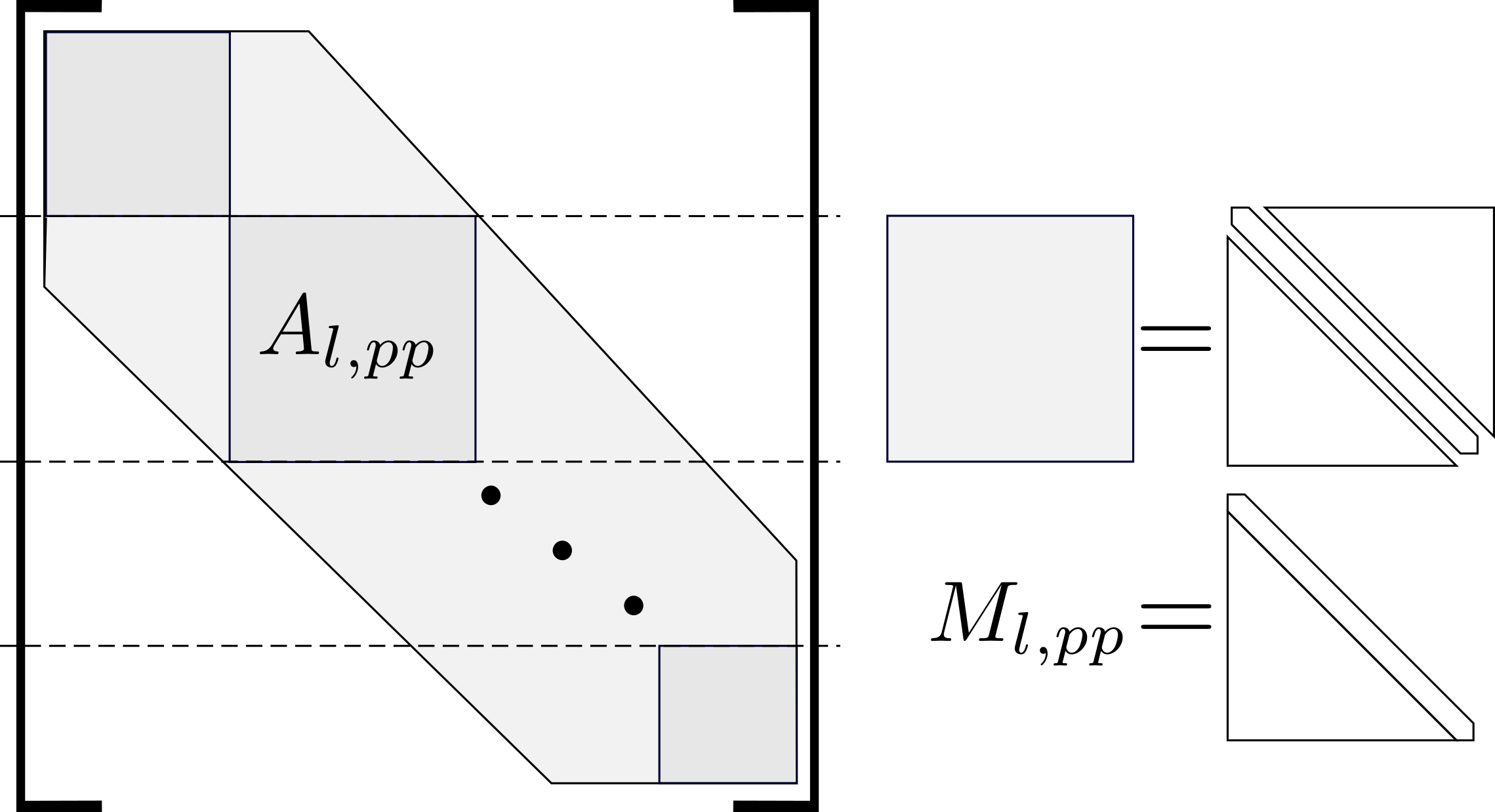}
    \caption{Depiction of the structure of the hybrid forward/backward Gauss-Seidel method on a general row-block parallel distribution of symmetric positive definite matrix $A$.}
    \label{fig:hfbgs}
\end{figure}

For the case at hand, we select each iteration matrix of the smoother sequence $\{M_l\}_{l=0}^{n_\ell-1}$ as the one representing four iterations $(\nu = 4)$ of the hybrid forward/backward Gauss-Seidel method. 
We consider having $A$ in a general row-block parallel distribution over 
$n_p$ processes, i.e., $A$ is divided into $n_p$ blocks of size $n_b \times 
n$, and we call $A_{pp}$ the corresponding diagonal block of~$A$. We then 
decompose each block $A_{pp}$ as 
$A_{pp}=L_{pp}+D_{pp}+L_{pp}^T$, where $D_{pp}=\operatorname{diag}(A_{pp})$, 
$L_{pp}$ is strictly lower triangular. To enforce symmetry 
in~\eqref{eq:multiplicative_composition}, we select $M_{l,pp}$ as the block 
diagonal matrices (Figure~\ref{fig:hfbgs})
\[
M_{l} = \operatorname{blockdiag}(M_{l,pp})_{p = 1}^{n_p/n_b}, \qquad M_{l,pp} = \omega \left( L_{l,pp}+D_{l,pp} \right),  \quad l=0,\ldots,n_\ell,
\]
where $\omega$ is a damping parameter. The overall procedure thus amounts essentially  to  using four sweeps of the damped block-Jacobi method on the 
matrix of the level while solving the blocks with the forward, respectively 
backward, Gauss-Seidel method.

\begin{figure}[htbp]
    \centering
    \subfloat[Adjacency matrix\label{fig:matching_example:adjacency}]{\includegraphics[width=0.3\columnwidth]{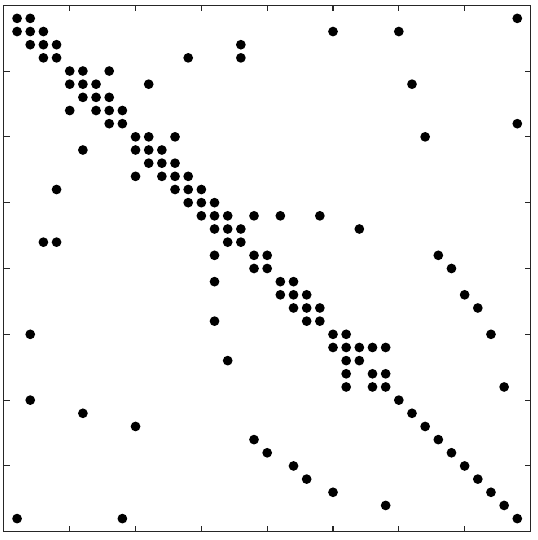}}
    \subfloat[Original graph\label{fig:matching_example:graph}]{\includegraphics[width=0.3\columnwidth]{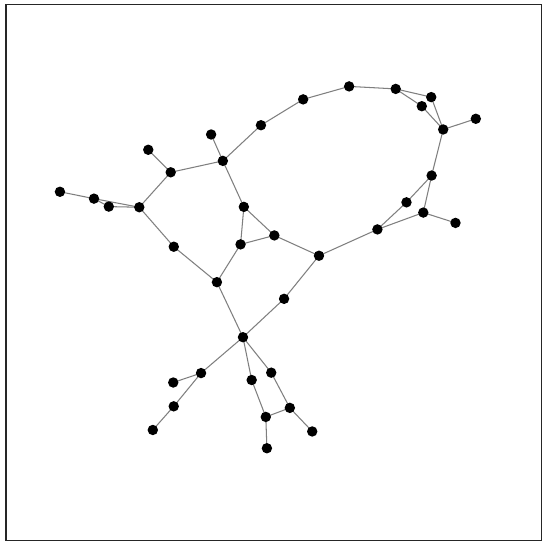}}
    \subfloat[Matched edges\label{fig:matching_example:matchededges}]{\includegraphics[width=0.3\columnwidth]{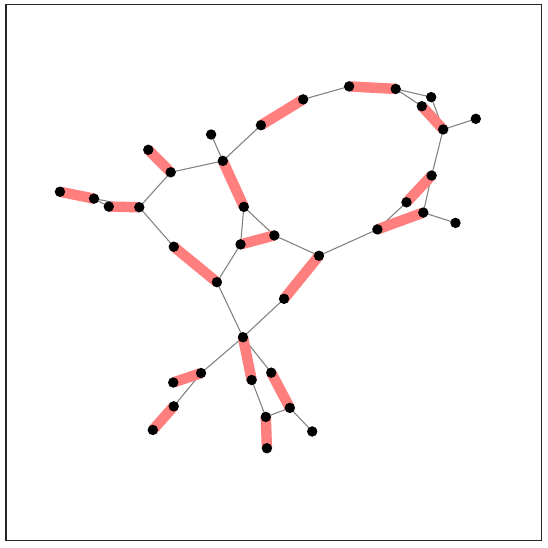}}
    \caption{Matching of the graph \texttt{bcspwr01} from the Harwell-Boeing collection. The matched nodes in the graph are highlighted by a bold red edge.}
    \label{fig:matching_example}
\end{figure}

To build the prolongation (and thus the restriction) matrices, we employ the 
{\em coarsening based on compatible weighted matching} strategy; a full 
account of the derivation and detailed theoretical analysis may be found 
in~\cite{DV2013,DFV2018,dambra2020b}.
This is a recursive procedure that starts from the adjacency graph $G = 
(V,E)$ associated with the sparse matrix $A$; this is the graph in which the 
vertex set $V$ consists of either the row or column indices of 
$A$ and the edge set $E$ corresponds to the indices pairs 
$(i,j)$ of the nonzero entries in $A$. The method works by 
constructing a \emph{matching} $\mathcal{M}$ in 
the graph $G$ to obtain a partition into subgraphs. We recall 
that a graph matching is  a subset of the graph's edges such 
that no two of them are incident on the same vertex. 
Specifically, we consider more than a purely 
topological matching by taking into account  the \emph{weights} of the 
edges, i.e., the values of the entries of the matrix $A$. In 
the  first step, we associate an edge weight matrix $C$, 
computed from the entries $a_{i,j}$ in $A$ and an arbitrary vector $\mathbf{w}$; then, we compute an 
approximate maximum product matching of the whole graph to obtain the 
aggregates defining the coarse spaces.  We define $C = (c_{i,j})_{i,j}$ as
\begin{equation}
\label{eq:weights}
c_{i,j} = 1 - \frac{2 a_{i,j} w_i w_j}{a_{i,i} w_i^2 + a_{j,j}w_j^2};
\end{equation}
then, $\mathcal{M}$ is an {\em approximate maximum product matching} of $G$ with edge weight matrix $C$, i.e., 
\begin{equation}
    \label{eq:approximatematch}
    \mathcal{M} \approx \arg \max_{\mathcal{M}'} \prod_{ (i,j) \in \mathcal{M}'} c_{i,j}.
\end{equation} 
The aggregates are then the subsets of indices $\{\mathcal{G}_p \}_{p=1}^{\lvert\mathcal{M}\rvert}$ of the whole index set $\mathcal{I}$ of $A$ made of pairs of indices matched by the algorithm, where  we denote with $\lvert\mathcal{M}\rvert$ the cardinality of the graph matching~$\mathcal{M}$. In other terms, we have obtained the decomposition
\[
\mathcal{I} = \{1,\ldots,n\} = \bigcup_{p=1}^{n_{\mathcal{M}}} \mathcal{G}_p, \quad \mathcal{G}_p \cap \mathcal{G}_r = \emptyset \text{ if } p\neq r;
\]
see, e.g., Figure~\ref{fig:matching_example} in which the matching of a test graph is computed{--in more detail, Figure~\ref{fig:matching_example:adjacency} has a black dot corresponding to a non-zero element of the adjacency matrix; Figure~\ref{fig:matching_example:graph} shows the corresponding graph obtained from it; while Figure~\ref{fig:matching_example:matchededges} highlights the aggregated nodes, i.e., the $\mathcal{G}_p$ sets.}
In most cases, we will end up with a sub-optimal matching, i.e., not all vertices will be endpoints of matched edges; thus, we usually have \emph{unmatched} vertices. To each unmatched vertex, we associate an aggregate $G_s$ that is a singleton, and we denote with $n_{\mathcal{S}}$ the total number of singletons. The main computational cost of this phase is represented by the computation of the \emph{approximate graph matching} on a graph that is distributed across thousands of processors. The parallel coarsening implemented in \texttt{AMG4PSBLAS} uses the \texttt{MatchBox-P} software library~\cite{CDGHP2011}; this implements a distributed parallel algorithm for the computation of \emph{half-approximate maximum weight matching} with complexity $\mathcal{O}(|E| \Delta)$, where $|E|$ is the cardinality of the graph edge set and $\Delta$ is the maximum vertex degree, i.e., the maximum number of edges incident on any given node of the graph. The  procedure guarantees a solution that is at least half of the optimal weight, i.e., the approximation in~\eqref{eq:approximatematch} holds within $1/2$ of the optimum. The message aggregation and overlapping between communication and computation employed by this strategy reduces the impact of the data communication on parallel efficiency; we refer the reader to~\cite{dambra2020b} for a complete set of experiments showcasing this feature. Finally, to build the prolongator matrices, the last ingredients we need are the vectors $\mathbf{w}_e$ identifying for each edge $e_{i \mapsto j} \in \mathcal{M}$  the orthonormal projection of $\mathbf{w}$ on the non-singleton aggregate $G_p$. For the sake of the explanation, we  consider an  ordering of the indices in which we move all the unknowns corresponding to unmatched vertices at the bottom\footnote{This ordering is for explanatory purposes only, and is not actually enforced in practice.}, and thus define a {\em tentative prolongator}
\begin{equation}\label{eq:prolongator}
\hat{P} = \left (
\begin{array}{cc}
\tilde{P} & 0\\
0 & W
\end{array}
\right ) \in \mathbb{R}^{n \times n_c}, 
\end{equation}
where:
\[
\tilde{P}=\operatorname{blockdiag}(w_{e_1}, \ldots, w_{e_{n_{\mathcal{M}}}} ),
\]
$W=\operatorname{diag}(w_s/|w_s|), \; s=1, \ldots, n_{\mathcal{S}}$, corresponds to unmatched vertices. The resulting number of coarse variables is then given by $n_c=n_{\mathcal{M}}+n_{\mathcal{S}}$. The matrix $\hat{P}$ we have just built is a piecewise constant interpolation operator whose range includes, by construction, the vector $\mathbf{w}$. The actual prolongator $P$ is then obtained from $\hat{P}$ as $P=(I-\omega D^{-1}A)\hat{P}$, where $D=\operatorname{diag}(A)$ and $\omega = 1/\|D^{-1}A\|_{\infty} \approx 1/\rho(D^{-1}A)$, with $\rho(D^{-1}A)$ the spectral radius of $D^{-1}A$. Indeed, the $P$ we have built is an instance of \emph{smoothed aggregation}. Please observe that the procedure we have described produces, at best, a halving of the size of the system at each new level of the hierarchy. Given the size of the systems we are interested in, this may be unsatisfactory since the number of levels in the hierarchy and thus the operational cost needed to cross it, would be too large. Fortunately, it is rather easy to overcome this issue: to obtain aggregates of size greater than two, we just have to collect them together by multiplying the corresponding prolongators (restrictors). This permits us to select the desired size of the aggregates (2, 4, 8, and so on) as an input parameter of the method.

To conclude the description of the preconditioners, we need to 
specify the choice for the coarsest solver. While using a 
direct solver at the coarsest  level is the easiest way to 
ensure that the coarsest grid is resolved to the needed 
tolerance, such an approach for an AMG method running on many 
thousands of parallel cores can be very expensive. If the 
matching strategy has worked satisfactorily, the 
coarsest-level matrix will tend to have both a small global 
size and a small number of rows per core: in this case the 
cost of data communication will dominate the local arithmetic  
computations causing a deterioration of the method efficiency. We use here a dual  strategy: on the one hand,  we employ a \emph{distributed coarsest solver} 
running on all the parallel cores, whilst on the other, we 
limit the maximum  size of  the coarsest-level matrix to 
around 200 unknowns per core. 
Specifically,  we use  the Flexible Conjugate Gradient (FCG) method with a block-Jacobi  preconditioner on which we solve 
approximately the blocks by an 
incomplete LU  factorization with one level of fill-in, 
ILU(1), the stopping criterion is  based on the reduction of the relative residual of $3$ orders 
of magnitude or a maximum number of iterations equal to $30$.

To have a comparison with the preconditioner just discussed, 
we also consider the same construction but with a different 
aggregation procedure: 
the decoupled version of the classic smoothed aggregation of 
\citet{MR1393006}. This is an aggregation option that was 
already available in previous versions of the 
library~\cite{MR2738234,Buttari2007223}, and was  already 
successfully used  in CFD 
applications~\cite{Aprovitola201115,Aprovitola20152688}. The 
basic idea is to build a coarse set of indices by grouping 
unknowns into disjoint subsets (the aggregates) by using an 
affinity measure and defining a simple tentative prolongator 
whose range contains the so-called near null 
space of the matrix of the given level, i.e., a sample of the eigenvector corresponding to the smallest eigenvalue. The strategy is implemented in an 
embarrassingly parallel fashion, i.e., each processor produces aggregates by only looking at local unknowns, i.e. the aggregation  is performed in a 
decoupled fashion, in contrast to the previous matching procedure that instead crosses the boundary of the single process.

Table~\ref{tab:preconditioners} summarizes the different preconditioners we have discussed here and that are used in the experiments of Section~\ref{comparison}.

\begin{table}[htbp]
    \centering\small
    \setlength{\tabcolsep}{0.2em}
    \begin{tabular}{lccc}
    \toprule
    Pre-smoother      & \multicolumn{3}{c}{$4$ iterations of hybrid forward Gauss-Seidel}   \\
    Post-smoother     & \multicolumn{3}{c}{$4$ iterations of hybrid backward Gauss-Seidel}  \\
    Coarsest solver   & \multicolumn{3}{c}{FCG preconditioned by block-Jacobi with ILU(1) block solvers}\\
    Cycle             & \multicolumn{3}{c}{V-cycle} \\
    Aggregation       & \multicolumn{2}{c}{Coupled smoothed based on matching} & Decoupled classic  \\
                      & $|\mathcal{G}| \leq 8$ & $|\mathcal{G}| \leq 16$ & smoothed\\
    \midrule
    Label & \emph{MLVSMATCH3} & \emph{MLVSMATCH4} & \emph{MLVSBM} \\
    \bottomrule
    \end{tabular}
    \caption{Summary of the described preconditioners, the labels are used to describe the results in Section~\ref{comparison}. }
    \label{tab:preconditioners}
\end{table}

\begin{remark}
The \texttt{AMG4PSBLAS} library provides interfaces to some widely used parallel direct solvers, such as \texttt{SuperLU}~\cite{L2005} and \texttt{MUMPS}~\cite{ADL2000}. Thus, we could have used any of those within the damped block-Jacobi method, either on the smoother or on the coarsest solvers. For what concerns the smoothers, it has been observed in the literature~\cite{BFKY2011,dambra2020b} that the combination with the Gauss-Seidel method delivers better smoothing properties for the overall method. In the coarsest solver case, the size of the local matrices is small enough to not usually show a significant performance increase when using a direct solver. We also stress that the preconditioner described in this section depends only on native \texttt{PSCToolkit} code, i.e., the user does not have to install optional third-party libraries to use it.
\end{remark}

\section{Interfacing Alya to \texttt{PSCToolkit}}
\label{interface}

The Alya code is organized in a modular way, and its architecture is split into \emph{modules}, \emph{kernel}, and \emph{services}, which can be separately compiled and linked. Each module represents a physical model, i.e., a set of partial differential equations which can interact for running a multi-physics simulation in a time-splitting approach, while Alya’s kernel implements the functionalities for dealing with the discretization mesh, the solvers and the I/O functionalities. As already mentioned, the governing equations of a physical model are discretized in space by using FE methods and all the functionalities to assemble the global stiffness matrix and right-hand-side (RHS) of the corresponding set of equations, including boundary conditions and material properties are the responsibility of the module. Instead, all the functionalities needed to solve the algebraic linear systems are implemented in the kernel. Some work on data structures and distribution of matrices and RHS was necessary to interface Alya with libraries from \texttt{PSCToolkit}, as described in the following.

\begin{figure*}[htb]
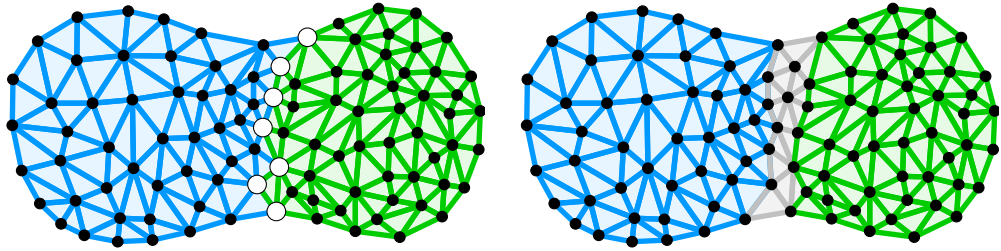

\centering
\subfloat[In white, interface nodes.\label{fig:meshleft}]{\input{mesh1.tikz}}
\subfloat[In white, halo elements.\label{fig:meshright}]{\input{mesh2.tikz}}

\caption{Mesh partitioning into (\ref{fig:meshleft}) disjoint sets of nodes, and (\ref{fig:meshright}) disjoint sets of elements. \label{meshdistribution}}
\end{figure*}

Alya uses the compressed sparse row matrix scheme for the internal representation of sparse matrices. This scheme is supported by \texttt{PSCToolkit}, so no significant difficulty was met from this perspective.
The main difficulty in the interfacing process was how the data, i.e. the discretization mesh and the corresponding unknowns, are distributed among the parallel processes and the way the related sparse matrix rows and RHS are locally assembled.
The Alya code is based on a domain decomposition where the discretization mesh is partitioned into disjoint subsets of elements/nodes, referred to as subdomains. Then, each subdomain is assigned to a parallel process that carries out all the geometrical and algebraic operations corresponding to that part of the domain and the associated unknowns. The interface elements/nodes on the boundary between two subdomains are assigned to one of the subdomains (see Figure~\ref{meshdistribution}).

The sparse matrices expressing the linear couplings among the unknowns are distributed in such a way that each parallel process holds the entries associated with the couplings generated on its subdomain. Two different options are possible for sparse matrix distribution: the partial row format and the full row format~\cite{alyaintechopen2017}, respectively. In the full row format, if a mesh element/node and the corresponding unknown belong to a process, all row entries related to that unknown are stored by that process. In the partial row format, the row of a matrix corresponding to an unknown is not full and needs contributions from unknowns belonging to different processes.%
Alya uses a partial row format for storing the matrix.

The libraries from \texttt{PSCToolkit} build the preconditioners and apply the Krylov methods on the assumption of a full row format;  nevertheless, support for partial row format was added to the libraries' pre-processing stage so that the interfacing can be as transparent as possible. 
The pre-processing support implies the retrieval of remote information for those matrix contributions that correspond to elements on the boundary; the data communication is split between the discovery of the needed entries (which needs only be executed when the discretization mesh changes) and the actual retrieval of the matrix entries, which must happen at any time step where the matrix coefficients and/or vector entries may be rebuilt, prior to an invocation of  the solvers. When the topology of the mesh does not change, and there is only an update in the coefficients, {it is also possible to reuse the same preconditioner; this may be full reuse of the overall matrices hierarchy, or partial reuse, employing the same prolongators/restrictors to rebuild the AMG hierarchy and smoothers.} 

{{We developed a software module in Alya's kernel for declaration, allocation, and initialization of the library's data structures as well as for using solvers and preconditioners.}
\texttt{PSBLAS} makes available some of the widely used iterative
methods based on Krylov projection methods through a single
interface to a driver routine, while preconditioners for \texttt{PSBLAS} Krylov solvers are available through the \texttt{AMG4PSBLAS} package. 
The main functionalities for selecting and building the chosen preconditioner are the responsibility of the software module included in the Alya's kernel, while the functionalities for applying it within the \texttt{PSBLAS} Krylov solver are completely transparent to the Alya code and are the responsibility of the library.

\section{The Bolund test case}
\label{bolund}

\begin{figure*}[htb]
\centering
\includegraphics[width=0.95\columnwidth]{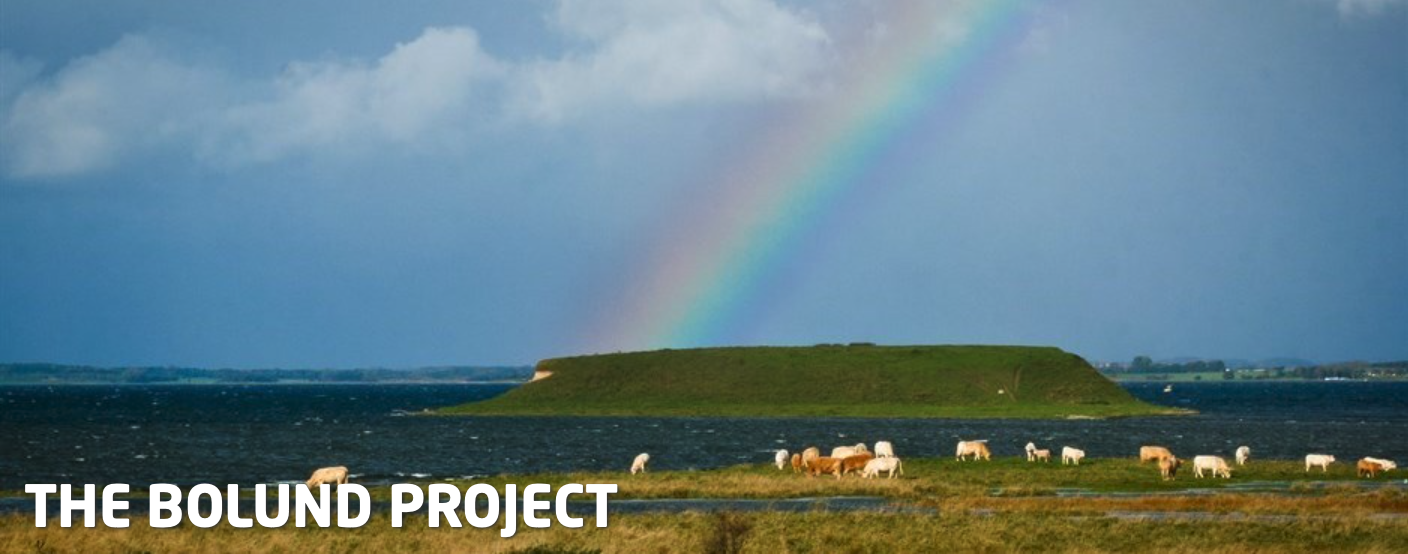}
\caption{Photograph of the Bolund hill~\cite{BolundWP}.\label{Bolund}}
\end{figure*}

Our main aim was to test the libraries for systems stemming from fluid dynamics simulation of incompressible flow arising in the study of wind-farm efficiency.
The test case is based on the Bolund experiment, a classical benchmark for microscale atmospheric flow models
over complex terrain~\cite{BSBMR2011,Bechmann2011}. {An incompressible flow treatment is used because the Mach number, i.e., the ratio of the speed of the flow to the speed of sound, is much smaller than $0.3$.}
The test case is based on a small (\SI{12}{\metre}) isolated steep hill at Roskilde Fjord in Denmark having a significantly steep escarpment in the main wind direction and uniformly covered by grass so that the resulting flow is not influenced by individual roughness elements. This is considered the ideal benchmark for the validation of neutral flow models and, hence a most relevant scenario for the analysis of software modeling for wind energy. Though relatively small, its geometrical shape induces complex 3D flow. Bolund was equipped with several measurement masts with conventional meteorological instruments and remote sensing Lidars to obtain detailed information of mean wind, wind shear, turbulence intensities, \emph{etc.} A publicly available database for evaluating currently available flow models and methodologies for turbine siting in complex terrain regarding wind resources and loads is available at~\cite{BolundWP}.

\begin{figure*}[htb]
\centering
\includegraphics[width=0.95\columnwidth]{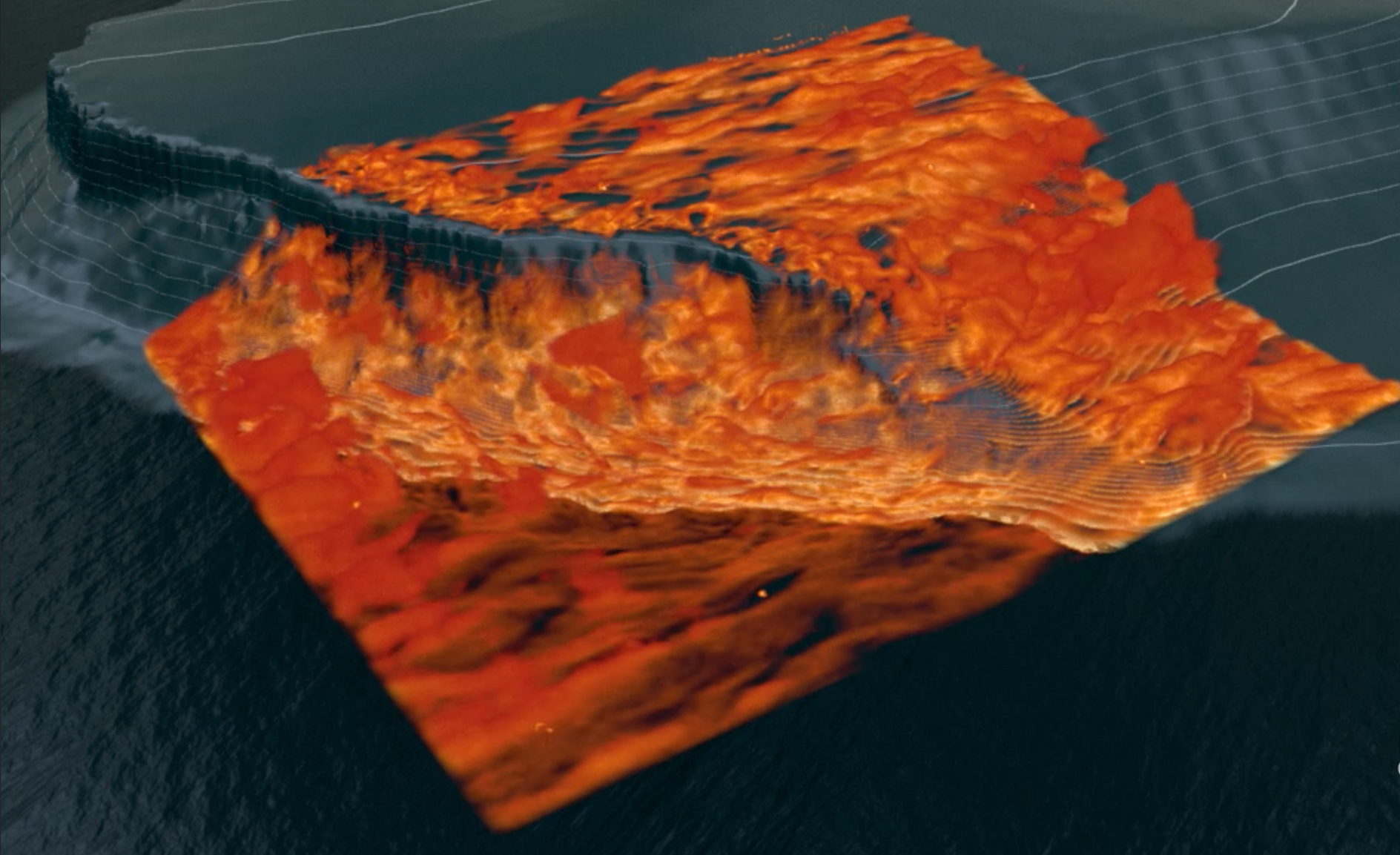}
\caption{Volume rendering of the velocity over Bolund obtained with Alya.\label{Bolund_alya}}
\end{figure*}

We discretize the incompressible Navier-Stokes equation~\eqref{eq:navierstokes} as described in Section~\ref{Alya Description}. At each time step of the LES procedure, we solved the SPD linear systems arising from the pressure equation employing the preconditioned flexible version of the CG method (FCG) method by \texttt{PSBLAS}. Starting from an initial guess for pressure from the previous time step, we stopped linear iterations when the Euclidean norm of the relative residual was no larger than $TOL=10^{-3}$. The Reynolds number based on the friction velocity for this test case is approximately $RE_\tau = U h / \nu \approx 10^7$ with $U = \SI{10}{\metre\per\second}$. As discussed in~\cite[Section~2.1]{BSBMR2011}, we can neglect Coriolis force in the horizontal direction and use the formulation~\eqref{eq:navierstokes} since the Rossby number $R_O = 667 \gg 1$.

The next Section~\ref{results} details the scalability result obtained for this test case with the new solvers and preconditioners from \texttt{PSCToolkit} described in Section~\ref{sec:PSCToolkit}.

\section{Parallel Performance Results}
\label{results}

In the following we discuss the results of experiments run on two of the most powerful European supercomputers. The first set of experiments aimed to analyze the behavior of different AMG preconditioners available from \texttt{AMG4PSBLAS} and run on the Marenostrum-4 supercomputer up to $12288$ CPU cores.
Marenostrum-4 is composed of $3456$ nodes with $2$ Intel Xeon Platinum $8160$ CPUs with $24$ cores per CPU. It is ranked 121\textsuperscript{th} in the November 2023 TOP500 list\footnote{Available at \href{https://www.top500.org}{www.top500.org}.}, with more than $10$ petaflops of peak performance and is operated by the Barcelona Supercomputer Center.
The simulations have been performed with the Alya code interfaced to \texttt{PSBLAS} (3.7.0.1) and \texttt{AMG4PSBLAS} (1.0), built with GNU compilers 7.2.
The second set of experiments aimed to reach very large scales and run by using only one of the most promising preconditioners by \texttt{AMG4PSBLAS} on the Juwels supercomputer, up to $23551$ CPU cores. Juwels is composed of $2271$ compute nodes with 2 Intel Xeon Platinum $8168$ CPUs, of $24$ cores each. It is ranked 127\textsuperscript{th} in the June 2022 TOP500 list, with more than $9$ petaflops of peak performance, and is operated by the J\"ulich Supercomputer Center. The simulations have been performed with the Alya code interfaced to the same versions of the solvers libraries mentioned above, built with GNU compilers 10.3.

\subsection{Comparison of AMG Preconditioners}
\label{comparison}

In this section, we discuss results obtained on Marenostrum-4 and compare the behavior of FCG coupled to the preconditioners described in Section~\ref{sec:amgprecond} and summarized in Table~\ref{tab:preconditioners}.
We run both strong scalability analysis for unstructured meshes of tetrahedra of three fixed sizes as well as weak scalability analysis, obtained by fixing different mesh sizes per core and linearly increasing both mesh size and the number of cores. 
A general row-block matrix distribution based on the Metis 4.0 mesh partitioner~\cite{MR1639073} was applied for the parallel runs.

\subsubsection{Strong scalability}
\label{strongscaling}

We first focus on strong scalability results obtained on the Bolund experiment for three fixed size meshes (small, medium and large) including $n_1=5570786 \approx 6\times 10^6$, $n_2=43619693 \approx 4.4\times 10^7$ and $n_3=345276325 \approx 0.35\times 10^9$ dofs, respectively.
{Three different configurations of the number of cores, obtained by doubling each time the number of MPI cores with respect to the minimum number of cores (nodes) needed to run at full load, were employed for the three different mesh sizes: from $min_p=48$ to $max_p=192$ cores in the case of the small mesh, from $min_p=384$ to $max_p=1536$ cores for the medium mesh, and finally from $min_p=3072$ to $max_p=12288$ cores for the large mesh}. We analyze the parallel efficiency and convergence behavior of the linear solvers for $20$ time steps after a pre-processing phase so that we focus on the solvers' behavior in the simulation of a fully developed flow. Note that in the Alya code a master-slave approach is employed, where the master process is not involved in the parallel computations.

In Figures~\ref{fig:strong1}-\ref{fig:strong1bis}, we report a comparison of the different methods in terms of the total number of iterations of the linear solvers and of the solve time per iteration (in seconds), respectively. Note that in the figures we also have results obtained with a version of Deflated CG ({\em AlyaDefCG}), available from the original Alya code.
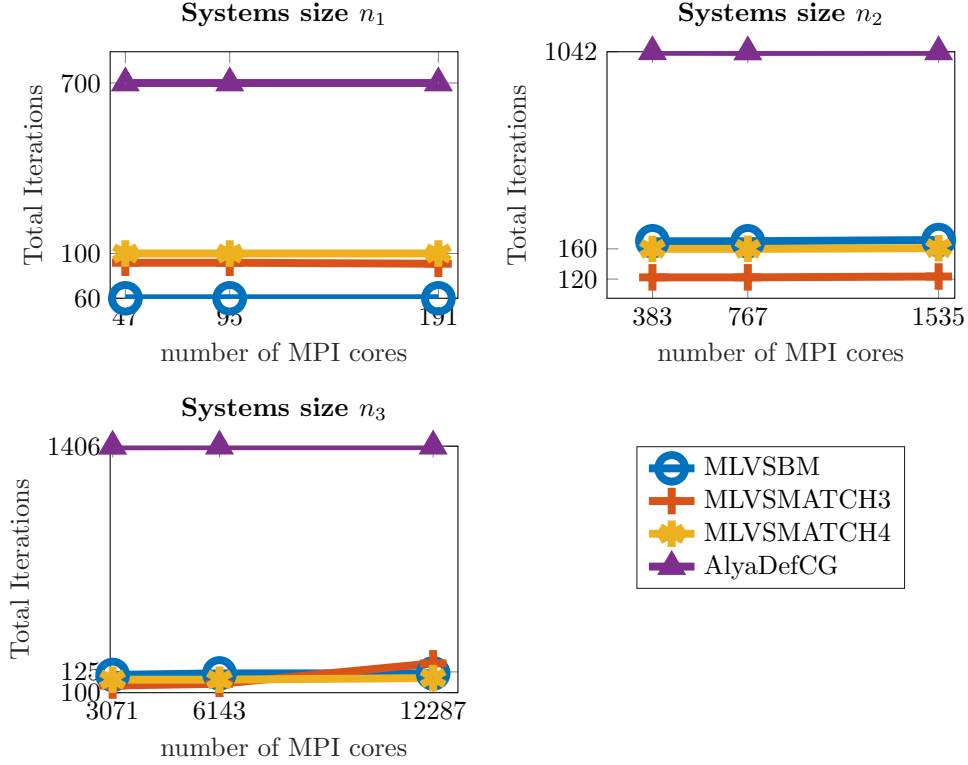
\begin{figure}[htbp]
\centering
\input{strongiterationmarenostrum.tikz}

\caption{Strong scalability: total iteration number of the linear solvers\label{fig:strong1}}
\end{figure}
\begin{figure}[htbp]
\centering
\input{strongtimeperiterationmarenostrum.tikz}

\caption{Strong scalability: time per iteration of the linear solvers\label{fig:strong1bis}}
\end{figure}
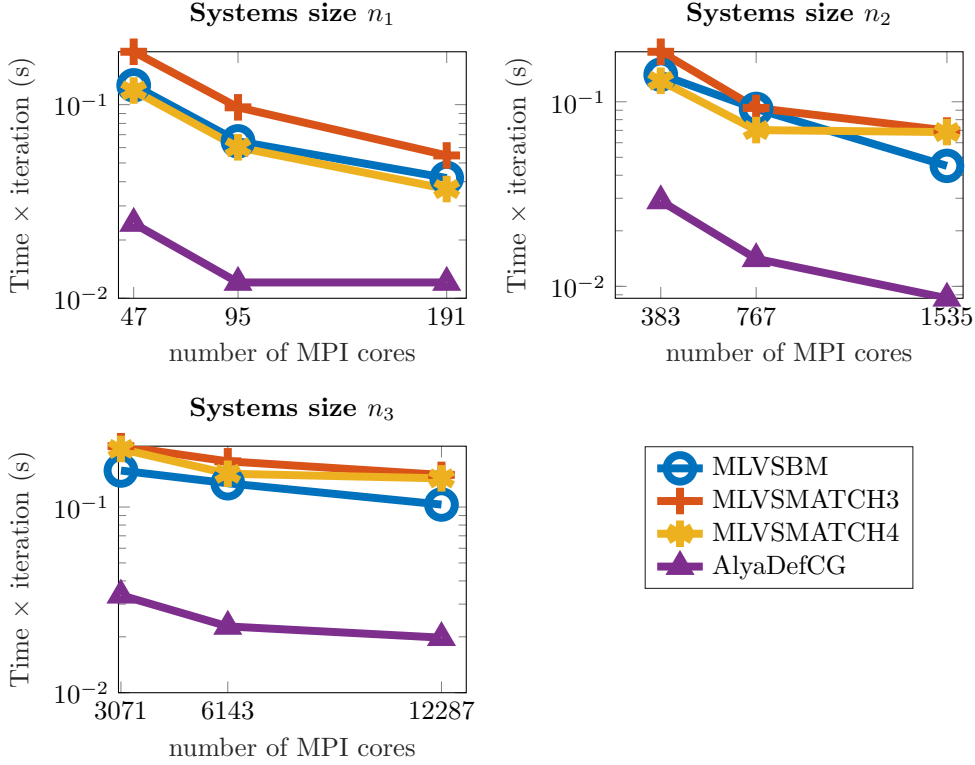

We can observe that the total number of linear iterations is much smaller than that with the original {\em AlyaDefCG}, for all three meshes, when 
\texttt{AMG4PSBLAS} multilevel preconditioners are applied. For the small mesh, the minimum number of linear iterations is obtained by {\em MLVSBM} 
which shows a fixed number of $60$ iterations for all core counts, while {\em MLVSMATCH3} requires $90$ iterations for all core counts except on $192$ cores, where $1$ less iteration was needed, and {\em MLVSMATCH4} 
requires $100$ iterations; in this case, the original {\em AlyaDefCG} requires $700$ iterations for all core counts.

In the case of the medium mesh, we observe a larger number of iterations of the solvers employing \texttt{AMG4PSBLAS} preconditioners with respect to the large mesh. We have a minimum number of iterations with {\em MLVSMATCH3} ranging from $122$ to $123$ for all number of cores, while {\em MLVSMATCH4} requires a range from $160$ to $161$ iterations and {\em MLVSBM} requires a range from $172$ to $174$ iterations. The original {\em AlyaDefCG} requires a number of iterations ranging from $1040$ to $1042$ for the medium mesh.

In the case of the large mesh, the number of iterations required by {\em MLVSMATCH3} ranges between $108$ on $3072$ cores and $137$ on $12288$ 
cores, while {\em MLVSMATCH4} requires a more stable number of iterations ranging from $115$ to $117$; a similarly stable behavior is observed for 
{\em MLVSBM} which requires a number of iterations ranging from $121$ to $123$. 
{\em AlyaDefCG} requires a number of iterations ranging from $1404$ to $1406$ for the large mesh.

The oscillations in the number of iterations seem to be mostly
dependent on the data partitioning obtained by Metis,  which 
in turn, appears  to have a larger impact on the {\em MLVSMATCH3} preconditioner in the case 
of the  large mesh. A deeper analysis of the impact of the 
data partitioner on the solver behavior, albeit interesting, 
is out of the scope of our current work and would require a 
significant amount of computing resources. 

In all cases, the time needed per iteration decreases for an 
increasing number of cores and, as expected, it is larger for 
the AMG preconditioners, where the cost for the preconditioner 
application at each FCG iteration is larger than that of {\em 
AlyaDefCG}. Depending on mesh size and number 
of cores, the AMG preconditioners show very similar behavior, although {\em MLVSBM} always requires a smaller time per iteration for the large mesh and for the medium mesh when $1536$ cores are used.

In Figures~\ref{fig:strong2}-\ref{fig:strong2bis}, we can see the total solve time spent in the linear solvers and the resulting speedup for the 
preconditioners. Here, we define speedup as the ratio $Sp=T_{min_p}/T_p$, where $T_{min_p}$ is the total time for solving linear systems when the 
minimum number of total cores, per each problem size, is involved in the simulation, and $T_p$ is the total time spent in linear solvers for all the 
increasing number of cores used for the specified mesh size.

We observe that the AMG preconditioners from 
\texttt{AMG4PSBLAS} generally achieve shorter execution times 
than the original {\em AlyaDefCG}; 
indeed, the expected longer time per iteration is more than 
compensated by the large reduction in the number of iterations 
especially for the small and 
large mesh. In good agreement with the behavior in terms of 
iterations and time per iteration, we observe that {\em 
MLVSBM} generally shows the shortest execution time for the 
small mesh, especially for small number of cores, while for the 
medium and large mesh, {\em MLVSMATCH3} and {\em MLVSMATCH4} 
show some better or comparable behaviour with respect to {\em MLVSBM}. The best speedups are generally obtained, except for 
the small mesh,  by the original {\em AlyaDefCG}, while in the case of AMG preconditioners, the very good convergence behavior and solve time on the smallest number of 
cores limit the speedup for the increasing number of cores. For the \texttt{AMG4PSBLAS} preconditioners, speedups are in good agreement with the total solve times, showing that {\em MLVSMATCH3} and {\em MLVSMATCH4} are generally better or comparable with respect to {\em MLVSBM} for all meshes when the small and medium number of cores are used, while {\em MLVSBM} is better for medium and large mesh when the largest number of cores is used.
\begin{figure}[htb]
\input{strongtotalsolvetimemarenostrum.tikz}
\caption{Strong scalability: total solve time of the linear solvers\label{fig:strong2}}
\end{figure}
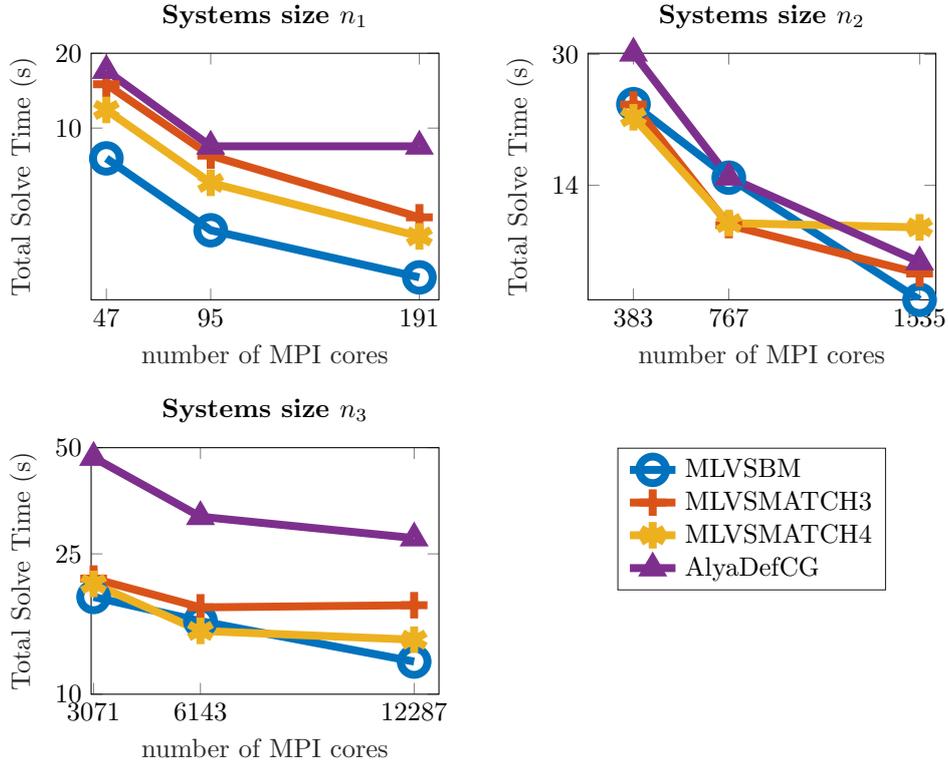
\begin{figure}[htb]
\centering
\input{strongspeedupofsolvemarenostrum.tikz}
\caption{Strong scalability: speedup of the linear solvers. {We note that ideal values for speedups in all three configurations are $1$, $2$ and $4$, respectively.}\label{fig:strong2bis}}
\end{figure}
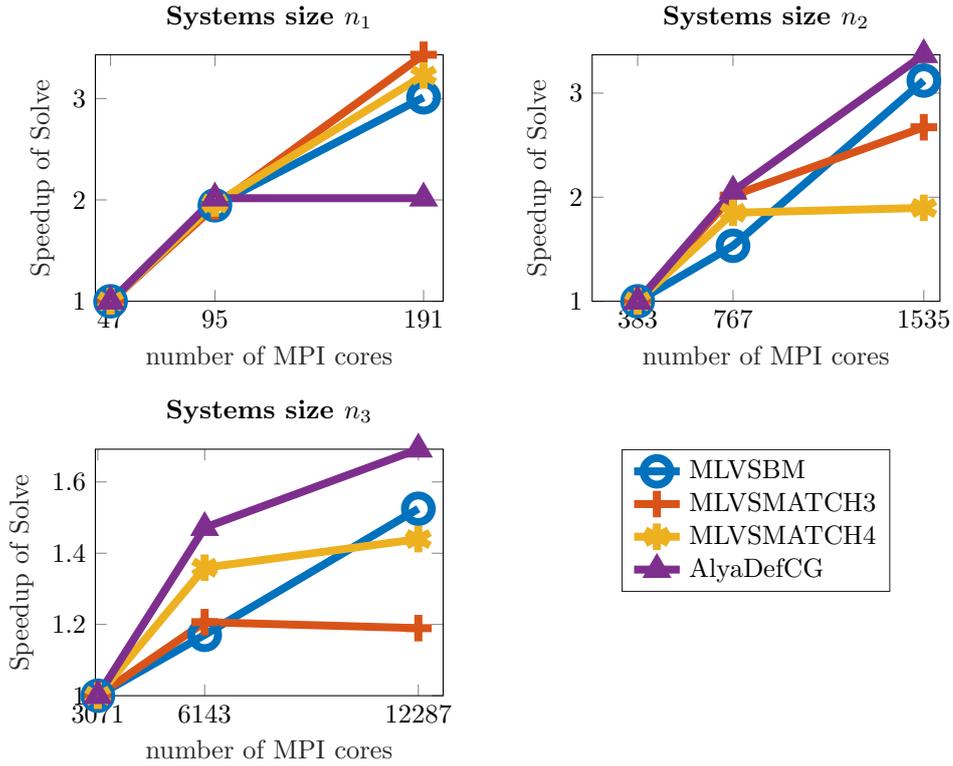

In conclusion, the selected solvers from the \texttt{PSCToolkit} generally outperform the original Alya solver for the employed test case, and the
choice of the better preconditioner from {\em AMG4PBLAS} depends on target mesh size and number of employed parallel cores. This appears as an advantage for Alya's users that having available a large set of parallel preconditioners through the interface to \texttt{PSCToolkit}, can
select the best one for their specific aims.

\subsubsection{Weak scalability}
\label{weakscaling}

In this section, we analyze the weak scalability of the \texttt{AMG4PSBLAS} 
preconditioners, i.e., we observe the solvers looking at their behavior when we fix the mesh size per core and increase the number of cores.

We considered the same test case and the three meshes of the previous section in the three possible configurations of computational cores, from $48$ up to $3072$, from $96$ up to $6144$ and from $192$ to $12288$. The different configurations of cores correspond to three different (decreasing) mesh sizes per core equal to $\text{nxcore}_1=1.1e5$, $\text{nxcore}_2=5.9e4$, and $\text{nxcore}_3=2.9e4$, respectively. Note that the medium and the large mesh correspond to scaling factors of $8$ and $64$, respectively, with respect to the small mesh; therefore in the same way, we scaled the number of cores for our weak scalability analysis. 

We can limit our analysis to observe the average number of linear iterations of the different 
employed preconditioners per each time step in the various simulations and to 
analyze execution times and scaled speedup for solve.
In Figure~\ref{fig:weak1}, we report the average number of iterations for each time step. We can observe a general increase, ranging from $35$ to $70$ 
for an increasing number of cores when the original {\em AlyaDefCG} is employed.
On the other hand, when AMG preconditioners from \texttt{AMG4PSBLAS} coupled with FCG by \texttt{PSBLAS} 
are applied, we observe a constant average number of iterations equal to $5$ for {\em MLVSMATCH4} both
for the small and the large mesh, independently of the number of cores, while {\em MLVSBM} requires $3$ iterations for the small mesh and $6$ for the large mesh. {\em MLVSMACTH3} ranges from $4$ to $6$ iterations on the small mesh and the large mesh, respectively. In the case of medium mesh, in agreement with what was observed for the strong scalability analysis, all the preconditioners require a larger average number of iterations, which is $8$ for {\em MLVSBM} and {\em MLVSMATCH4}, and $6$ for {\em MLVSMATCH3}.
This behavior indicates a very promising algorithmic scalability of {\em MLVSMATCH4}.
\begin{figure}[htb]
\centering
\input{weakavgiterationmarenostrum.tikz}
\caption{Weak scalability: average number of linear iterations per time step. $\text{nxcore}_1$ dofs per core (\ref{fig:weak1-top}), $\text{nxcore}_2$ dofs per core (\ref{fig:weak1-middle}), $\text{nxcore}_3$ dofs per core (\ref{fig:weak1-bottom})\label{fig:weak1}}
\end{figure}
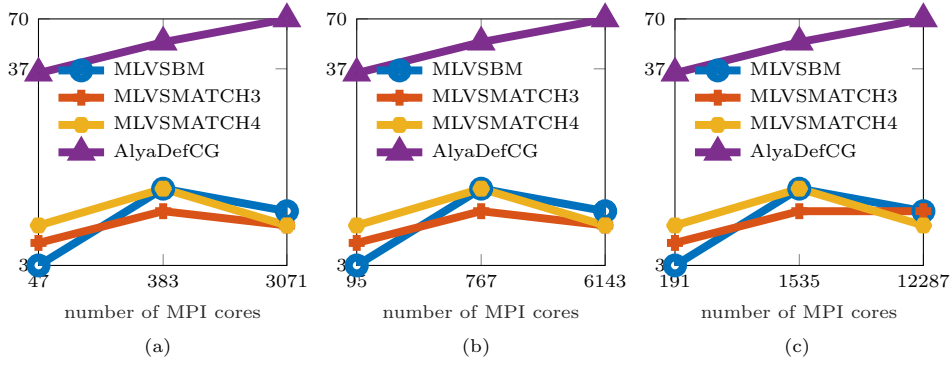
In Figures~\ref{fig:weak2quater}-\ref{fig:weak2tris}, we can see the total solve time and the corresponding scaled speedup. We can observe that, as expected from the previous sections, all preconditioners from \texttt{AMG4PSBLAS} generally lead to a smaller increase ratio in the solve times with respect to the original {\em AlyaDefCG}, when the mesh size goes from the small to the large one. In more detail, we observe that, for all mesh sizes per core, smaller increase ratios in the execution time are generally obtained with {\em MLVSMATCH3} and {\em MLVSMATCH4}. This is better observed by looking at the scaled speedup. It is defined as $\text{scalfactor}\times T_{min_p}/T_p$, where $\text{scalfactor}=1,8,64$, for the three increasing number of cores, $T_{min_p}$ is the total time for solving linear systems when the minimum number of total cores is involved in the simulation, per each mesh size per core, and $T_p$ is the total time spent in linear solvers for all the increasing number of cores used for the specified mesh size per core. We observe that the best values are obtained with the {\em MLVSMATCH3} and {\em MLVSMATCH4} preconditioners when $\text{nxcore}_1$ and $\text{nxcore}_2$ dofs per core are used. In detail, for $\text{nxcore}_1$ dofs per core, {\em MLVSMATCH3} reaches the best value of about $71 \%$ of scaled efficiency on $3072$ cores and about $44 \%$ of scaled efficiency on $6144$ core when $\text{nxcore}_2$ dofs per core are employed. This shows that the scalability of {\em MLVSMATCH3} and {\em MLVSMATCH4} are very promising in facing the exascale challenge, especially when the resources are used at their best in terms of node memory capacity and bandwidth. On the other hand, in the case of $\text{nxcore}_3$ dofs per core (\ref{weak2tris-bottom}), the scaled speedup of {\em AlyaDefCG} is better; this is essentially due to the very large solve time spent by this solver on $192$ cores.
\begin{figure}[htb]
\input{weaktotalsolvetimemarenostrum.tikz}
\caption{Weak scalability: total solve time~(s) of the linear solvers. $\text{nxcore}_1$ dofs per core (\ref{fig:weak2quater-top}), $\text{nxcore}_2$ dofs per core (\ref{fig:weak2quater-middle}), $\text{nxcore}_3$ dofs per core (\ref{fig:weak2quater-bottom})\label{fig:weak2quater}}
\end{figure}
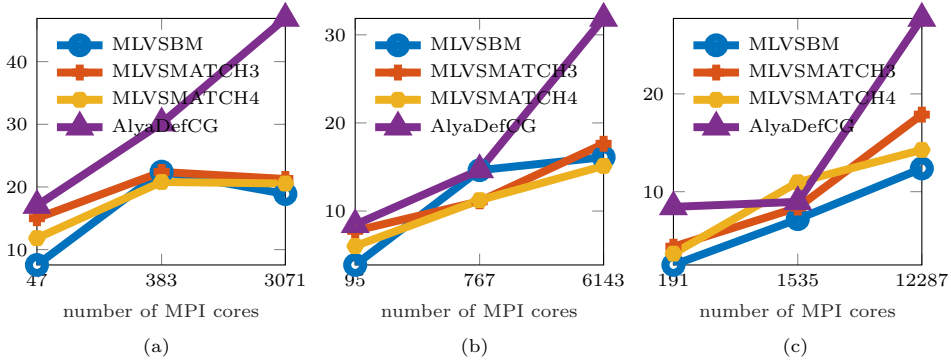
\begin{figure}[htb]
\centering
\input{weakspeedupofsolvemarenostrum.tikz}

\caption{Weak scalability: the scaled speedup of the linear solvers. $\text{nxcore}_1$ dofs per core (\ref{weak2tris-top}), $\text{nxcore}_2$ dofs per core (\ref{weak2tris-middle}), $\text{nxcore}_3$ dofs per core (\ref{weak2tris-bottom})\label{fig:weak2tris}}
\end{figure}
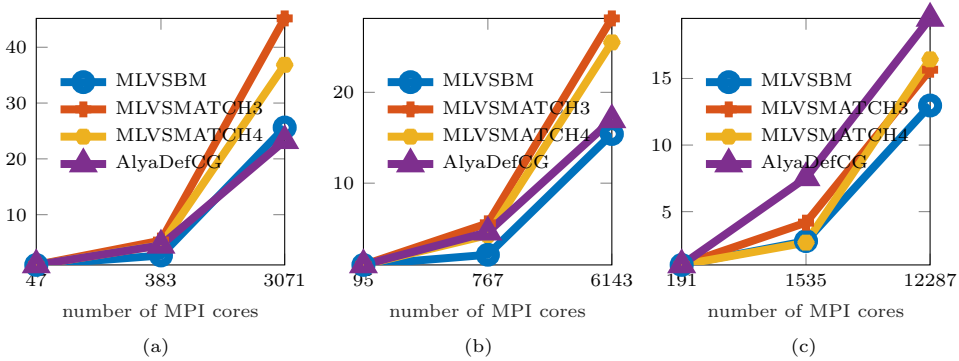
\subsection{Results at extreme scales}
\label{extremescales}

In this section, we discuss some results obtained on the Juwels supercomputer by increasing the number of dofs till to $n_4 \approx 2.9\times 10^9$. We limit our analysis to the weak scalability results of one of the most promising solvers in \texttt{PSCToolkit}. Indeed, due to the limited access to the Juwels resources and taking into account the above preconditioners comparison, we only run experiments by using the {\em MLVSMATCH4} preconditioner.
A general row-block data distribution based on a parallel geometric partitioning
using Space Filling Curve (SFC)~\cite{Borrelletal2018} was applied for these experiments.
As in the previous experiments, we analyze the parallel efficiency and convergence behavior of the linear solver for $20$ time steps after a pre-processing phase so that we focus on the solver behavior in the simulation of a fully developed flow for all the meshes but the largest one, where we were not able to skip the transient phase due to long simulation time. In this last case, we considered a total number of time steps equal to $1379$ and analyzed solver performance in the last $20$ time steps. Note that increasing mesh size imposes a decrease in time step due to stability constraints of the explicit time discretization that is preferred for LES simulations. Therefore, the total simulated time depends on the mesh size.
Furthermore, to reduce observed operating oscillations associated with the full node runs, we used only a total of $46$ cores per node.

As already mentioned, we analyze the weak scalability of the solvers;
we considered a mesh size per core equal to $\text{nxcore}_1$ and used a scaling factor of $8$ for going up to the largest mesh size; therefore in the same way, we scaled the number of cores for our weak scalability analysis.
We can limit our analysis to observing the average number of linear iterations of the solver per each time step and analyzing execution times and scaled speedup for the solve phase. We compare the results obtained by using the \texttt{PSCToolkit}'s solver against Alya's Conjugate Gradient solver (hereby {\em AlyaCG}). Observe that in these experiments, we also tried to use the Deflated CG implemented in Alya, but it does not work for the two larger test cases, and {\em AlyaCG} appears better in the case of smaller size meshes.
In Figure~\ref{fig:weak1bis}, we report the average number of iterations per each time step.
We can observe a general increase, ranging from $133$ to $331$ for an increasing number of cores, but on $368$ cores where $95$ average iteration count is obtained, when the original {\em AlyaCG} is employed, while very good algorithmic scalability, with an average number of linear iterations per each time step ranging from $4$ to $6$, when the \texttt{PSCToolkit}'s solver is applied.
\begin{figure}[htb]
\centering
\input{weakavgiterationjuwels.tikz}
\caption{Weak scalability: average number of linear iterations per time step. Systems size from $n_1$ to $n_4$. \label{fig:weak1bis}}
\end{figure}
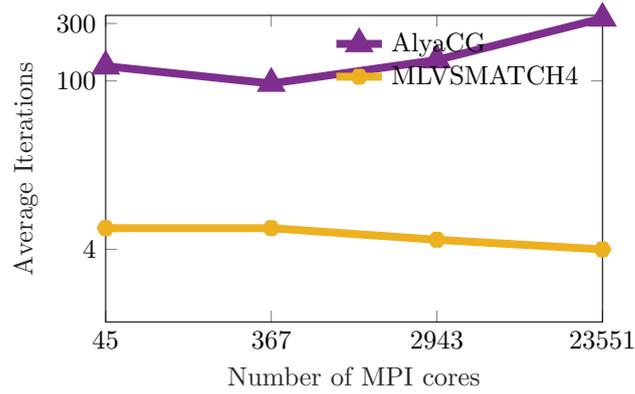
In Figures~\ref{fig:weak2}-\ref{fig:weak2bis}, we can see the total solve time and the corresponding scaled speedup.
We can observe that the good algorithmic scalability of {\em MLVSMATCH4} leads to an almost flat execution time for solving 
when the first three meshes are employed, while a decrease is 
observed for the simulation carried out with the largest mesh, 
depending on a smaller average number of iterations per time 
step. On the contrary, the original {\em AlyaCG} generally 
shows a huge increase for increasing number of cores 
and mesh size, but in the second one, where a decrease in the 
average number of iterations per time step is observed.
Then we look at the scaled speedup, defined as 
$\text{scalfactor} \times T_{45}/T_p$, where 
$\text{scalfactor}=1,8,64,512$, for increasing number of 
cores, $T_{45}$ is the total time for solving linear systems 
when $45$ cores are involved in the simulation, and $T_p$ is 
the total time spent in linear solvers for all the increasing 
number of cores. We observe that for the two larger meshes, 
{\em MLVSMATCH4} has a super-linear scaled speed-up 
{of about $71$ (up from the ideal speedup of  
$64$) and $640$ (up from the ideal speedup of $512$), 
respectively,} showing that its very good algorithmic 
scalability is coupled with excellent implementation 
scalability of all the basic computational kernels. This 
scalability is very promising in facing the exascale challenge.
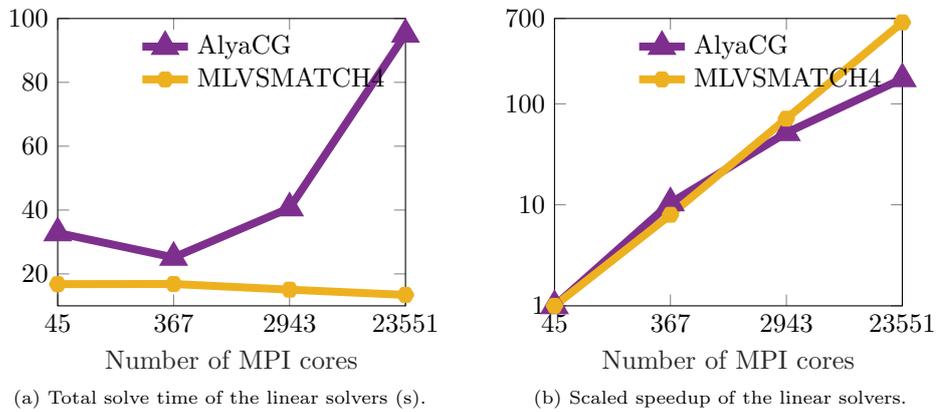
\begin{figure}[htb]
\subfloat[Total solve time of the linear solvers~(s). \label{fig:weak2}]{\input{weaktotalsolvetimejuwels.tikz}}\hfil
\subfloat[Scaled speedup of the linear solvers.\label{fig:weak2bis}]{\input{weakspeedupofsolvejuwels.tikz}}
\caption{Weak scalability: systems size from $n_1$ to $n_4$.}
\end{figure}

\section{Conclusions}\label{sec:conclusions}

In this paper, we presented our work on improving the linear 
solver  capabilities of a large-scale CFD code by
interfacing it with a software framework, including new and 
state-of-the-art  algebraic linear solvers, specifically 
designed to  exploit the very large   
potential of current petascale supercomputers and aimed at the 
early  exascale supercomputers. Our activities were carried 
out in the context of the European Center of Excellence for 
Energy applications, where one of the lighthouse codes was the 
Alya code, developed at the Barcelona 
Supercomputing Center (BSC) and applied to wind flow studies 
for renewable  production. However, this work has a wider 
impact, and confirms the benefits of using third-party 
software libraries developed by specialists, in 
complex, multi-component and multi-physics simulation codes. 

From Alya’s perspective, the most significant achievement has 
been obtaining  excellent algorithmic scalability thanks to 
multigrid preconditioners, as  shown in the weak scalability 
studies. This allows us to solve much bigger problems 
efficiently. Our first objective for the future is to test the 
GPU  version of \texttt{PSCToolkit}. During EoCoE, we have 
significantly  optimized the FE assembly on GPUs, making it 
four times more energy  efficient than the CPU version; 
integrating  a competitive linear algebra 
GPU package is now the next priority. 
After that, having the entire workflow for 
incompressible flow problems on GPUs should be relatively straightforward.
We expect to have a 
much higher number of unknowns for problems 
running for MPI process when GPU accelerators are exploited. 
Therefore, strong scalability should be much less critical.
While we have focused on a wind energy problem in this work, 
we wish to test  the solver in other incompressible flow 
problems in the future. Moreover, since the solver is fully 
interfaced with Alya, it will be interesting to 
test the suitability of \texttt{PSCToolkit} for other 
problems, such as solid mechanics or heat transfer.

\bmhead{Acknowledgments} 

We thank the two anonymous reviewers whose suggestions helped improve and clarify this manuscript. P.D., F.D. and S.F are members of the INdAM Research Group GNCS.

\section{Declarations}
 
\subsection*{Ethical Approval} N/A  
 
\subsection*{Competing interests} N/A  
 
\subsection*{Authors' contributions}  
The first two authors are involved in the development of the Alya code and worked on the integration of the PSCToolkit into Alya. Herbert Owen did main contribution in setup and running of experiments. The last three authors are the developers of PSCToolkit and experts in Linear Algebra. They have equally contributed in setup of suitable linear solvers, analyzing results and writing the paper.     
 
\subsection*{Funding} 
The research received funding and PRACE grants for supercomputers access from Horizon 2020 Project  ``Energy oriented Centre of Excellence: toward exascale for energy'' (EoCoE--II), Project ID: 824158.
 
\subsection*{Availability of data and materials}  
PSCToolkit is available at https://psctoolkit.github.io/
The data that support the findings of this study are available from the corresponding author upon reasonable request.

\bibliography{cfdbibliography.bib}

\end{document}

%% file: strongiterationmarenostrum.tikz
\definecolor{mycolor1}{rgb}{0.00000,0.44700,0.74100}%
\definecolor{mycolor2}{rgb}{0.85000,0.32500,0.09800}%
\definecolor{mycolor3}{rgb}{0.92900,0.69400,0.12500}%
\definecolor{mycolor4}{rgb}{0.49400,0.18400,0.55600}%
\begin{tikzpicture}

\begin{axis}[%
width=0.35\columnwidth,
height=0.25\columnwidth,
at={(0in,0in)},
scale only axis,
xmin=40,
xmax=200,
xtick={ 47,  95, 191},
xlabel style={font=\color{white!15!black}},
xlabel={number of MPI cores},
ymode=log,
ymin=60,
ymax=1000,
yminorticks=true,
ylabel style={font=\color{white!15!black}},
ylabel={Total Iterations},
ytick={60,100,700},
yticklabels={60,100,700},
axis background/.style={fill=white},
title style={font=\bfseries},
title={Systems size $n_1$},
legend style={at={(-0.03,1)}, anchor=north east, legend cell align=left, align=left, draw=white!15!black}
]
\addplot [color=mycolor1, line width=3.0pt, mark size=5.0pt, mark=o, mark options={solid, mycolor1}]
  table[row sep=crcr]{%
47	60\\
95	60\\
191	60\\
};
%\addlegendentry{MLVSBM}

\addplot [color=mycolor2, line width=3.0pt, mark size=5.0pt, mark=+, mark options={solid, mycolor2}]
  table[row sep=crcr]{%
47	90\\
95	90\\
191	89\\
};
%\addlegendentry{MLVSMATCH3}

\addplot [color=mycolor3, line width=3.0pt, mark size=5.0pt, mark=asterisk, mark options={solid, mycolor3}]
  table[row sep=crcr]{%
47	100\\
95	100\\
191	100\\
};
%\addlegendentry{MLVSMATCH4}

\addplot [color=mycolor4, line width=3.0pt, mark size=3.3pt, mark=triangle, mark options={solid, mycolor4}]
  table[row sep=crcr]{%
47	700\\
95	700\\
191	700\\
};
%\addlegendentry{AlyaDefCG}

\end{axis}

\begin{axis}[%
width=0.35\columnwidth,
height=0.25\columnwidth,
at={(0.5\columnwidth,0in)},
scale only axis,
xmin=200,
xmax=1600,
xtick={ 383,  767, 1535},
xticklabels={ 383,  767, 1535},
xlabel style={font=\color{white!15!black}},
xlabel={number of MPI cores},
ymode=log,
ymin=100,
ymax=1042,
ytick={120,160,1042},
yticklabels={120,160,1042},
yminorticks=true,
ylabel style={font=\color{white!15!black}},
ylabel={Total Iterations},
axis background/.style={fill=white},
title style={font=\bfseries},
title={Systems size $n_2$},
legend style={at={(-0.03,1)}, anchor=north east, legend cell align=left, align=left, draw=white!15!black}
]
\addplot [color=mycolor1, line width=3.0pt, mark size=5.0pt, mark=o, mark options={solid, mycolor1}]
  table[row sep=crcr]{%
383	172\\
767	172\\
1535	174\\
};
%\addlegendentry{MLVSBM}

\addplot [color=mycolor2, line width=3.0pt, mark size=5.0pt, mark=+, mark options={solid, mycolor2}]
  table[row sep=crcr]{%
383	122\\
767	122\\
1535	123\\
};
%\addlegendentry{MLVSMATCH3}

\addplot [color=mycolor3, line width=3.0pt, mark size=5.0pt, mark=asterisk, mark options={solid, mycolor3}]
  table[row sep=crcr]{%
383	160\\
767	160\\
1535	161\\
};
%\addlegendentry{MLVSMATCH4}

\addplot [color=mycolor4, line width=3.0pt, mark size=3.3pt, mark=triangle, mark options={solid, mycolor4}]
  table[row sep=crcr]{%
383	1042\\
767	1040\\
1535	1041\\
};
%\addlegendentry{AlyaDefCG}

\end{axis}

\begin{axis}[%
width=0.35\columnwidth,
height=0.25\columnwidth,
at={(0in,-0.4\columnwidth)},
scale only axis,
xmin=3000,
xmax=13000,
xtick={ 3071,  6143, 12287},
xticklabels={ 3071,  6143, 12287},
scaled x ticks = false,
xlabel style={font=\color{white!15!black}},
xlabel={number of MPI cores},
ymode=log,
ymin=100,
ymax=1406,
yminorticks=true,
ylabel style={font=\color{white!15!black}},
ylabel={Total Iterations},
ytick={100,125,1406},
yticklabels={100,125,1406},
axis background/.style={fill=white},
title style={font=\bfseries},
title={Systems size $n_3$},
legend style={at={(0.8\columnwidth,1)}, anchor=north east, legend cell align=left, align=left, draw=white!15!black}
]
\addplot [color=mycolor1, line width=3.0pt, mark size=5.0pt, mark=o, mark options={solid, mycolor1}]
  table[row sep=crcr]{%
3071	121\\
6143	123\\
12287	123\\
};
\addlegendentry{MLVSBM}

\addplot [color=mycolor2, line width=3.0pt, mark size=5.0pt, mark=+, mark options={solid, mycolor2}]
  table[row sep=crcr]{%
3071	108\\
6143	110\\
12287	137\\
};
\addlegendentry{MLVSMATCH3}

\addplot [color=mycolor3, line width=3.0pt, mark size=5.0pt, mark=asterisk, mark options={solid, mycolor3}]
  table[row sep=crcr]{%
3071	115\\
6143	115\\
12287	117\\
};
\addlegendentry{MLVSMATCH4}

\addplot [color=mycolor4, line width=3.0pt, mark size=3.3pt, mark=triangle, mark options={solid, mycolor4}]
  table[row sep=crcr]{%
3071	1404\\
6143	1406\\
12287	1406\\
};
\addlegendentry{AlyaDefCG}

\end{axis}
\end{tikzpicture}%

%% file: strongtimeperiterationmarenostrum.tikz
\definecolor{mycolor1}{rgb}{0.00000,0.44700,0.74100}%
\definecolor{mycolor2}{rgb}{0.85000,0.32500,0.09800}%
\definecolor{mycolor3}{rgb}{0.92900,0.69400,0.12500}%
\definecolor{mycolor4}{rgb}{0.49400,0.18400,0.55600}%
\begin{tikzpicture}

\begin{axis}[%
width=0.35\columnwidth,
height=0.25\columnwidth,
at={(0in,0in)},
scale only axis,
xmin=40,
xmax=200,
xtick={ 47,  95, 191},
xlabel style={font=\color{white!15!black}},
xlabel={number of MPI cores},
ymode=log,
ymin=0.01,
ymax=0.18765,
yminorticks=true,
ylabel style={font=\color{white!15!black}},
ylabel={Time $\times$ iteration (s)},
axis background/.style={fill=white},
title style={font=\bfseries},
title={Systems size $n_1$},
legend style={at={(-0.03,1)}, anchor=north east, legend cell align=left, align=left, draw=white!15!black}
]
\addplot [color=mycolor1, line width=3.0pt, mark size=5.0pt, mark=o, mark options={solid, mycolor1}]
  table[row sep=crcr]{%
47	0.1259\\
95	0.0646666666666667\\
191	0.0418\\
};
%\addlegendentry{MLVSBM}

\addplot [color=mycolor2, line width=3.0pt, mark size=5.0pt, mark=+, mark options={solid, mycolor2}]
  table[row sep=crcr]{%
47	0.18765\\
95	0.096875\\
191	0.054675\\
};
%\addlegendentry{MLVSMATCH3}

\addplot [color=mycolor3, line width=3.0pt, mark size=5.0pt, mark=asterisk, mark options={solid, mycolor3}]
  table[row sep=crcr]{%
47	0.11832\\
95	0.06022\\
191	0.03666\\
};
%\addlegendentry{MLVSMATCH4}

\addplot [color=mycolor4, line width=3.0pt, mark size=3.3pt, mark=triangle, mark options={solid, mycolor4}]
  table[row sep=crcr]{%
47	0.0243428571428571\\
95	0.0120685714285714\\
191	0.0120685714285714\\
};
%\addlegendentry{AlyaDefCG}

\end{axis}

\begin{axis}[%
width=0.35\columnwidth,
height=0.25\columnwidth,
at={(0.5\columnwidth,0in)},
scale only axis,
xmin=200,
xmax=1600,
xtick={ 383,  767, 1535},
xticklabels={ 383,  767, 1535},
xlabel style={font=\color{white!15!black}},
xlabel={number of MPI cores},
ymode=log,
ymin=0.00861730769230769,
ymax=0.18665,
yminorticks=true,
ylabel style={font=\color{white!15!black}},
ylabel={Time $\times$ iteration (s)},
axis background/.style={fill=white},
title style={font=\bfseries},
title={Systems size $n_2$},
legend style={at={(-0.03,1)}, anchor=north east, legend cell align=left, align=left, draw=white!15!black}
]
\addplot [color=mycolor1, line width=3.0pt, mark size=5.0pt, mark=o, mark options={solid, mycolor1}]
  table[row sep=crcr]{%
383	0.140275\\
767	0.0914875\\
1535	0.04495\\
};
%\addlegendentry{MLVSBM}

\addplot [color=mycolor2, line width=3.0pt, mark size=5.0pt, mark=+, mark options={solid, mycolor2}]
  table[row sep=crcr]{%
383	0.18665\\
767	0.0925833333333333\\
1535	0.06985\\
};
%\addlegendentry{MLVSMATCH3}

\addplot [color=mycolor3, line width=3.0pt, mark size=5.0pt, mark=asterisk, mark options={solid, mycolor3}]
  table[row sep=crcr]{%
383	0.13\\
767	0.07025\\
1535	0.0685125\\
};
%\addlegendentry{MLVSMATCH4}

\addplot [color=mycolor4, line width=3.0pt, mark size=3.3pt, mark=triangle, mark options={solid, mycolor4}]
  table[row sep=crcr]{%
383	0.0290307692307692\\
767	0.0141076923076923\\
1535	0.00861730769230769\\
};
%\addlegendentry{AlyaDefCG}

\end{axis}

\begin{axis}[%
width=0.35\columnwidth,
height=0.25\columnwidth,
at={(0in,-0.4\columnwidth)},
scale only axis,
xmin=3000,
xmax=13000,
xtick={ 3071,  6143, 12287},
xticklabels={ 3071,  6143, 12287},
scaled x ticks = false,
xlabel style={font=\color{white!15!black}},
xlabel={number of MPI cores},
ymode=log,
ymin=0.01,%0.0197957142857143,
ymax=0.21272,
yminorticks=true,
ylabel style={font=\color{white!15!black}},
ylabel={Time $\times$ iteration (s)},
axis background/.style={fill=white},
title style={font=\bfseries},
title={Systems size $n_3$},
legend style={at={(0.8\columnwidth,1)}, anchor=north east, legend cell align=left, align=left, draw=white!15!black}
]
\addplot [color=mycolor1, line width=3.0pt, mark size=5.0pt, mark=o, mark options={solid, mycolor1}]
  table[row sep=crcr]{%
3071	0.15725\\
6143	0.1343\\
12287	0.103116666666667\\
};
\addlegendentry{MLVSBM}

\addplot [color=mycolor2, line width=3.0pt, mark size=5.0pt, mark=+, mark options={solid, mycolor2}]
  table[row sep=crcr]{%
3071	0.21272\\
6143	0.17632\\
12287	0.149\\
};
\addlegendentry{MLVSMATCH3}

\addplot [color=mycolor3, line width=3.0pt, mark size=5.0pt, mark=asterisk, mark options={solid, mycolor3}]
  table[row sep=crcr]{%
3071	0.2055\\
6143	0.15118\\
12287	0.14274\\
};
\addlegendentry{MLVSMATCH4}

\addplot [color=mycolor4, line width=3.0pt, mark size=3.3pt, mark=triangle, mark options={solid, mycolor4}]
  table[row sep=crcr]{%
3071	0.0334957142857143\\
6143	0.0227571428571429\\
12287	0.0197957142857143\\
};
\addlegendentry{AlyaDefCG}

\end{axis}
\end{tikzpicture}%

%% file: strongtotalsolvetimemarenostrum.tikz
\definecolor{mycolor1}{rgb}{0.00000,0.44700,0.74100}%
\definecolor{mycolor2}{rgb}{0.85000,0.32500,0.09800}%
\definecolor{mycolor3}{rgb}{0.92900,0.69400,0.12500}%
\definecolor{mycolor4}{rgb}{0.49400,0.18400,0.55600}%
\begin{tikzpicture}

\begin{axis}[%
width=0.35\columnwidth,
height=0.25\columnwidth,
at={(0in,0in)},
scale only axis,
xmin=40,
xmax=200,
xtick={ 47,  95, 191},
xlabel style={font=\color{white!15!black}},
xlabel={number of MPI cores},
ymin=0,
ymax=20,
ylabel style={font=\color{white!15!black}},
ylabel={Total Solve Time (s)},
ymode=log,
ytick={2,10,20},
yticklabels={2,10,20},
axis background/.style={fill=white},
title style={font=\bfseries},
title={Systems size $n_1$},
legend style={at={(-0.03,1)}, anchor=north east, legend cell align=left, align=left, draw=white!15!black}
]
\addplot [color=mycolor1, line width=3.0pt, mark size=5.0pt, mark=o, mark options={solid, mycolor1}]
  table[row sep=crcr]{%
47	7.5543\\
95	3.8804\\
191	2.5079\\
};
%\addlegendentry{MLVSBM}

\addplot [color=mycolor2, line width=3.0pt, mark size=5.0pt, mark=+, mark options={solid, mycolor2}]
  table[row sep=crcr]{%
47	15.0121\\
95	7.7508\\
191	4.3741\\
};
%\addlegendentry{MLVSMATCH3}

\addplot [color=mycolor3, line width=3.0pt, mark size=5.0pt, mark=asterisk, mark options={solid, mycolor3}]
  table[row sep=crcr]{%
47	11.8324\\
95	6.0218\\
191	3.6659\\
};
%\addlegendentry{MLVSMATCH4}

\addplot [color=mycolor4, line width=3.0pt, mark size=3.3pt, mark=triangle, mark options={solid, mycolor4}]
  table[row sep=crcr]{%
47	17.0401\\
95	8.44837\\
191	8.44837\\
};
%\addlegendentry{AlyaDefCG}

\end{axis}

\begin{axis}[%
width=0.35\columnwidth,
height=0.25\columnwidth,
at={(0.5\columnwidth,0in)},
scale only axis,
xmin=200,
xmax=1600,
xtick={ 383,  767, 1535},
xticklabels={383,  767, 1535},
xlabel style={font=\color{white!15!black}},
xlabel={number of MPI cores},
ymin=7.19220018386841,
ymax=30.1910991668701,
ymode=log,
ytick={7,14,30},
yticklabels={7,14,30},
ylabel style={font=\color{white!15!black}},
ylabel={Total Solve Time (s)},
axis background/.style={fill=white},
title style={font=\bfseries},
title={Systems size $n_2$},
legend style={at={(-0.03,1)}, anchor=north east, legend cell align=left, align=left, draw=white!15!black}
]
\addplot [color=mycolor1, line width=3.0pt, mark size=5.0pt, mark=o, mark options={solid, mycolor1}]
  table[row sep=crcr]{%
383	22.4442\\
767	14.6382\\
1535	7.1922\\
};
%\addlegendentry{MLVSBM}

\addplot [color=mycolor2, line width=3.0pt, mark size=5.0pt, mark=+, mark options={solid, mycolor2}]
  table[row sep=crcr]{%
383	22.397\\
767	11.11\\
1535	8.3821\\
};
%\addlegendentry{MLVSMATCH3}

\addplot [color=mycolor3, line width=3.0pt, mark size=5.0pt, mark=asterisk, mark options={solid, mycolor3}]
  table[row sep=crcr]{%
383	20.7992\\
767	11.2393\\
1535	10.9611\\
};
%\addlegendentry{MLVSMATCH4}

\addplot [color=mycolor4, line width=3.0pt, mark size=3.3pt, mark=triangle, mark options={solid, mycolor4}]
  table[row sep=crcr]{%
383	30.1911\\
767	14.6718\\
1535	8.96248\\
};
%\addlegendentry{AlyaDefCG}

\end{axis}

\begin{axis}[%
width=0.35\columnwidth,
height=0.25\columnwidth,
at={(0in,-0.4\columnwidth)},
scale only axis,
xmin=3000,
xmax=13000,
xtick={ 3071,  6143, 12287},
xticklabels={ 3071,  6143, 12287},
scaled x ticks = false,
xlabel style={font=\color{white!15!black}},
xlabel={number of MPI cores},
ymin=10,
ymax=50,
ymode=log,
ytick={10,25,50},
yticklabels={10,25,50},
ylabel style={font=\color{white!15!black}},
ylabel={Total Solve Time (s)},
axis background/.style={fill=white},
title style={font=\bfseries},
title={Systems size $n_3$},
legend style={at={(0.8\columnwidth,1)}, anchor=north east, legend cell align=left, align=left, draw=white!15!black}
]
\addplot [color=mycolor1, line width=3.0pt, mark size=5.0pt, mark=o, mark options={solid, mycolor1}]
  table[row sep=crcr]{%
3071	18.8699\\
6143	16.1164\\
12287	12.375\\
};
\addlegendentry{MLVSBM}

\addplot [color=mycolor2, line width=3.0pt, mark size=5.0pt, mark=+, mark options={solid, mycolor2}]
  table[row sep=crcr]{%
3071	21.2721\\
6143	17.6321\\
12287	17.8802\\
};
\addlegendentry{MLVSMATCH3}

\addplot [color=mycolor3, line width=3.0pt, mark size=5.0pt, mark=asterisk, mark options={solid, mycolor3}]
  table[row sep=crcr]{%
3071	20.5496\\
6143	15.1177\\
12287	14.2742\\
};
\addlegendentry{MLVSMATCH4}

\addplot [color=mycolor4, line width=3.0pt, mark size=3.3pt, mark=triangle, mark options={solid, mycolor4}]
  table[row sep=crcr]{%
3071	46.8942\\
6143	31.8603\\
12287	27.7146\\
};
\addlegendentry{AlyaDefCG}

\end{axis}
\end{tikzpicture}%

%% file: strongspeedupofsolvemarenostrum.tikz
\definecolor{mycolor1}{rgb}{0.00000,0.44700,0.74100}%
\definecolor{mycolor2}{rgb}{0.85000,0.32500,0.09800}%
\definecolor{mycolor3}{rgb}{0.92900,0.69400,0.12500}%
\definecolor{mycolor4}{rgb}{0.49400,0.18400,0.55600}%
\begin{tikzpicture}

\begin{axis}[%
width=0.35\columnwidth,
height=0.25\columnwidth,
at={(0in,0in)},
scale only axis,
xmin=40,
xmax=200,
xtick={ 47,  95, 191},
xlabel style={font=\color{white!15!black}},
xlabel={number of MPI cores},
%ymode=log,
ymin=1,
ymax=3.43204316316499,
yminorticks=true,
ylabel style={font=\color{white!15!black}},
ylabel={Speedup of Solve},
axis background/.style={fill=white},
title style={font=\bfseries},
title={Systems size $n_1$},
legend style={at={(-0.03,1)}, anchor=north east, legend cell align=left, align=left, draw=white!15!black}
]
\addplot [color=mycolor1, line width=3.0pt, mark size=5.0pt, mark=o, mark options={solid, mycolor1}]
  table[row sep=crcr]{%
47	1\\
95	1.94678383671786\\
191	3.01220144343873\\
};
%\addlegendentry{MLVSBM}

\addplot [color=mycolor2, line width=3.0pt, mark size=5.0pt, mark=+, mark options={solid, mycolor2}]
  table[row sep=crcr]{%
47	1\\
95	1.9368452288796\\
191	3.43204316316499\\
};
%\addlegendentry{MLVSMATCH3}

\addplot [color=mycolor3, line width=3.0pt, mark size=5.0pt, mark=asterisk, mark options={solid, mycolor3}]
  table[row sep=crcr]{%
47	1\\
95	1.96492743033644\\
191	3.22769306309501\\
};
%\addlegendentry{MLVSMATCH4}

\addplot [color=mycolor4, line width=3.0pt, mark size=3.3pt, mark=triangle, mark options={solid, mycolor4}]
  table[row sep=crcr]{%
47	1\\
95	2.01696895377451\\
191	2.01696895377451\\
};
%\addlegendentry{AlyaDefCG}

\end{axis}

\begin{axis}[%
width=0.35\columnwidth,
height=0.25\columnwidth,
at={(0.5\columnwidth,0in)},
scale only axis,
xmin=200,
xmax=1600,
xtick={ 383,  767, 1535},
xticklabels={ 383,  767, 1535},
xlabel style={font=\color{white!15!black}},
xlabel={number of MPI cores},
%ymode=log,
ymin=1,
ymax=3.36861002758165,
yminorticks=true,
ylabel style={font=\color{white!15!black}},
ylabel={Speedup of Solve},
axis background/.style={fill=white},
title style={font=\bfseries},
title={Systems size $n_2$},
legend style={at={(-0.03,1)}, anchor=north east, legend cell align=left, align=left, draw=white!15!black}
]
\addplot [color=mycolor1, line width=3.0pt, mark size=5.0pt, mark=o, mark options={solid, mycolor1}]
  table[row sep=crcr]{%
383	1\\
767	1.53326228634668\\
1535	3.12063068324018\\
};
%\addlegendentry{MLVSBM}

\addplot [color=mycolor2, line width=3.0pt, mark size=5.0pt, mark=+, mark options={solid, mycolor2}]
  table[row sep=crcr]{%
383	1\\
767	2.01593159315932\\
1535	2.67200343589315\\
};
%\addlegendentry{MLVSMATCH3}

\addplot [color=mycolor3, line width=3.0pt, mark size=5.0pt, mark=asterisk, mark options={solid, mycolor3}]
  table[row sep=crcr]{%
383	1\\
767	1.85057788296424\\
1535	1.89754677906414\\
};
%\addlegendentry{MLVSMATCH4}

\addplot [color=mycolor4, line width=3.0pt, mark size=3.3pt, mark=triangle, mark options={solid, mycolor4}]
  table[row sep=crcr]{%
383	1\\
767	2.05776387355335\\
1535	3.36861002758165\\
};
%\addlegendentry{AlyaDefCG}

\end{axis}

\begin{axis}[%
width=0.35\columnwidth,
height=0.25\columnwidth,
at={(0in,-0.4\columnwidth)},
scale only axis,
xmin=3000,
xmax=13000,
xtick={ 3071,  6143, 12287},
xticklabels={ 3071,  6143, 12287},
scaled x ticks = false,
xlabel style={font=\color{white!15!black}},
xlabel={number of MPI cores},
%ymode=log,
ymin=1,
ymax=1.69203957480895,
yminorticks=true,
ylabel style={font=\color{white!15!black}},
ylabel={Speedup of Solve},
axis background/.style={fill=white},
title style={font=\bfseries},
title={Systems size $n_3$},
legend style={at={(0.8\columnwidth,1)}, anchor=north east, legend cell align=left, align=left, draw=white!15!black}
]
\addplot [color=mycolor1, line width=3.0pt, mark size=5.0pt, mark=o, mark options={solid, mycolor1}]
  table[row sep=crcr]{%
3071	1\\
6143	1.17085081035467\\
12287	1.5248404040404\\
};
\addlegendentry{MLVSBM}

\addplot [color=mycolor2, line width=3.0pt, mark size=5.0pt, mark=+, mark options={solid, mycolor2}]
  table[row sep=crcr]{%
3071	1\\
6143	1.20644166038078\\
12287	1.1897014574781\\
};
\addlegendentry{MLVSMATCH3}

\addplot [color=mycolor3, line width=3.0pt, mark size=5.0pt, mark=asterisk, mark options={solid, mycolor3}]
  table[row sep=crcr]{%
3071	1\\
6143	1.35930730203669\\
12287	1.439632343669\\
};
\addlegendentry{MLVSMATCH4}

\addplot [color=mycolor4, line width=3.0pt, mark size=3.3pt, mark=triangle, mark options={solid, mycolor4}]
  table[row sep=crcr]{%
3071	1\\
6143	1.47186937976102\\
12287	1.69203957480895\\
};
\addlegendentry{AlyaDefCG}

\end{axis}
\end{tikzpicture}%

%% file: weakavgiterationmarenostrum.tikz
\definecolor{mycolor1}{rgb}{0.00000,0.44700,0.74100}%
\definecolor{mycolor2}{rgb}{0.85000,0.32500,0.09800}%
\definecolor{mycolor3}{rgb}{0.92900,0.69400,0.12500}%
\definecolor{mycolor4}{rgb}{0.49400,0.18400,0.55600}%
\subfloat[\label{fig:weak1-top}]{\begin{tikzpicture}
\begin{axis}[%
width=0.25\columnwidth,
height=0.25\columnwidth,
at={(2.5in,7.1in)},
scale only axis,
xmin=47,
xmax=3071,
xtick={  47,  383, 3071},
xticklabels={  47,  383, 3071},
yticklabel style = {font=\footnotesize},
xticklabel style = {font=\footnotesize},
xmode=log,
xlabel style={font=\color{white!15!black}\footnotesize},
xlabel={number of MPI cores},
ymin=3,
ymax=70,
ymode=log,
ytick={3,37,70},
yticklabels={3,37,70},
ylabel style={font=\color{white!15!black}\footnotesize},
axis background/.style={fill=white},
title style={font=\small},
legend style={at={(0.03,0.87)}, anchor=north west, legend cell align=left, align=left, draw=none, fill=none, font=\footnotesize}
]
\addplot [color=mycolor1, line width=3.0pt, mark size=3.0pt, mark=o, mark options={solid, mycolor1}]
  table[row sep=crcr]{%
47	3\\
383	8\\
3071	6\\
};
\addlegendentry{MLVSBM}

\addplot [color=mycolor2, line width=3.0pt, mark size=3.0pt, mark=+, mark options={solid, mycolor2}]
  table[row sep=crcr]{%
47	4\\
383	6\\
3071	5\\
};
\addlegendentry{MLVSMATCH3}

\addplot [color=mycolor3, line width=3.0pt, mark size=3.0pt, mark=asterisk, mark options={solid, mycolor3}]
  table[row sep=crcr]{%
47	5\\
383	8\\
3071	5\\
};
\addlegendentry{MLVSMATCH4}

\addplot [color=mycolor4, line width=3.0pt, mark size=3.3pt, mark=triangle, mark options={solid, mycolor4}]
  table[row sep=crcr]{%
47	35\\
383	52\\
3071	70\\
};
\addlegendentry{AlyaDefCG}

\end{axis}
\end{tikzpicture}}
\subfloat[\label{fig:weak1-middle}]{\begin{tikzpicture}
\begin{axis}[%
width=0.25\columnwidth,
height=0.25\columnwidth,
at={(2.5in,4.101in)},
scale only axis,
xmin=95,
xmax=6143,
xtick={  95,  767, 6143},
xticklabels={  95,  767, 6143},
yticklabel style = {font=\footnotesize},
xticklabel style = {font=\footnotesize},
xlabel style={font=\color{white!15!black}\footnotesize},
xlabel={number of MPI cores},
ymin=3,
ymax=70,
xmode=log,
ymode=log,
ytick={3,37,70},
yticklabels={3,37,70},
ylabel style={font=\color{white!15!black}\footnotesize},
axis background/.style={fill=white},
title style={font=\small},
legend style={at={(0.03,0.87)}, anchor=north west, legend cell align=left, align=left, draw=none, fill=none, font=\footnotesize}
]
\addplot [color=mycolor1, line width=3.0pt, mark size=3.0pt, mark=o, mark options={solid, mycolor1}]
  table[row sep=crcr]{%
95	3\\
767	8\\
6143	6\\
};
\addlegendentry{MLVSBM}

\addplot [color=mycolor2, line width=3.0pt, mark size=3.0pt, mark=+, mark options={solid, mycolor2}]
  table[row sep=crcr]{%
95	4\\
767	6\\
6143	5\\
};
\addlegendentry{MLVSMATCH3}

\addplot [color=mycolor3, line width=3.0pt, mark size=3.0pt, mark=asterisk, mark options={solid, mycolor3}]
  table[row sep=crcr]{%
95	5\\
767	8\\
6143	5\\
};
\addlegendentry{MLVSMATCH4}

\addplot [color=mycolor4, line width=3.0pt, mark size=3.3pt, mark=triangle, mark options={solid, mycolor4}]
  table[row sep=crcr]{%
95	35\\
767	52\\
6143	70\\
};
\addlegendentry{AlyaDefCG}

\end{axis}
\end{tikzpicture}}
\subfloat[\label{fig:weak1-bottom}]{\begin{tikzpicture}
\begin{axis}[%
width=0.25\columnwidth,
height=0.25\columnwidth,
at={(2.5in,1.101in)},
scale only axis,
xmin=191,
xmax=12287,
xmode=log,
xtick={  191,  1535, 12287},
xticklabels={  191,  1535, 12287},
scaled x ticks = false,
yticklabel style = {font=\footnotesize},
xticklabel style = {font=\footnotesize},
xlabel style={font=\color{white!15!black}\footnotesize},
xlabel={number of MPI cores},
ymin=3,
ymax=70,
ymode=log,
ytick={3,37,70},
yticklabels={3,37,70},
ylabel style={font=\color{white!15!black}\footnotesize},
axis background/.style={fill=white},
title style={font=\small},
legend style={at={(0.03,0.87)}, anchor=north west, legend cell align=left, align=left, draw=none, fill=none, font=\footnotesize}
]
\addplot [color=mycolor1, line width=3.0pt, mark size=3.0pt, mark=o, mark options={solid, mycolor1}]
  table[row sep=crcr]{%
191	3\\
1535	8\\
12287	6\\
};
\addlegendentry{MLVSBM}

\addplot [color=mycolor2, line width=3.0pt, mark size=3.0pt, mark=+, mark options={solid, mycolor2}]
  table[row sep=crcr]{%
191	4\\
1535	6\\
12287	6\\
};
\addlegendentry{MLVSMATCH3}

\addplot [color=mycolor3, line width=3.0pt, mark size=3.0pt, mark=asterisk, mark options={solid, mycolor3}]
  table[row sep=crcr]{%
191	5\\
1535	8\\
12287	5\\
};
\addlegendentry{MLVSMATCH4}

\addplot [color=mycolor4, line width=3.0pt, mark size=3.3pt, mark=triangle, mark options={solid, mycolor4}]
  table[row sep=crcr]{%
191	35\\
1535	52\\
12287	70\\
};
\addlegendentry{AlyaDefCG}

\end{axis}
\end{tikzpicture}}%

%% file: weaktotalsolvetimemarenostrum.tikz
\definecolor{mycolor1}{rgb}{0.00000,0.44700,0.74100}%
\definecolor{mycolor2}{rgb}{0.85000,0.32500,0.09800}%
\definecolor{mycolor3}{rgb}{0.92900,0.69400,0.12500}%
\definecolor{mycolor4}{rgb}{0.49400,0.18400,0.55600}%
\subfloat[\label{fig:weak2quater-top}]{\begin{tikzpicture}
\begin{axis}[%
width=0.25\columnwidth,
height=0.25\columnwidth,
at={(2.486in,7.1in)},
scale only axis,
xmin=47,
xmax=3071,
xtick={  47,  383, 3071},
xticklabels={47,383,3071},
yticklabel style = {font=\footnotesize},
xticklabel style = {font=\footnotesize},
xmode=log,
xlabel style={font=\color{white!15!black}\footnotesize},
xlabel={number of MPI cores},
ymin=7.5543,
ymax=46.8942,
ylabel style={font=\color{white!15!black}\footnotesize},
axis background/.style={fill=white},
title style={font=\small},
legend style={at={(0.03,0.97)}, anchor=north west, legend cell align=left, align=left, draw=none, fill=none, font=\footnotesize}
]
\addplot [color=mycolor1, line width=3.0pt, mark size=3.0pt, mark=o, mark options={solid, mycolor1}]
  table[row sep=crcr]{%
47	7.5543\\
383	22.4442\\
3071	18.8699\\
};
\addlegendentry{MLVSBM}

\addplot [color=mycolor2, line width=3.0pt, mark size=3.0pt, mark=+, mark options={solid, mycolor2}]
  table[row sep=crcr]{%
47	15.0121\\
383	22.397\\
3071	21.2721\\
};
\addlegendentry{MLVSMATCH3}

\addplot [color=mycolor3, line width=3.0pt, mark size=3.0pt, mark=asterisk, mark options={solid, mycolor3}]
  table[row sep=crcr]{%
47	11.8324\\
383	20.7992\\
3071	20.5496\\
};
\addlegendentry{MLVSMATCH4}

\addplot [color=mycolor4, line width=3.0pt, mark size=3.3pt, mark=triangle, mark options={solid, mycolor4}]
  table[row sep=crcr]{%
47	17.0401\\
383	30.1911\\
3071	46.8942\\
};
\addlegendentry{AlyaDefCG}

\end{axis}
\end{tikzpicture}}
\subfloat[\label{fig:weak2quater-middle}]{\begin{tikzpicture}
\begin{axis}[%
width=0.25\columnwidth,
height=0.25\columnwidth,
at={(2.486in,4.101in)},
scale only axis,
xmin=95,
xmax=6143,
xtick={  95,  767, 6143},
xticklabels={  95,  767, 6143},
yticklabel style = {font=\footnotesize},
xticklabel style = {font=\footnotesize},
xmode=log,
xlabel style={font=\color{white!15!black}\footnotesize},
xlabel={number of MPI cores},
ymin=3.8804,
ymax=31.8603,
ylabel style={font=\color{white!15!black}\footnotesize},
axis background/.style={fill=white},
title style={font=\small},
legend style={at={(0.03,0.97)}, anchor=north west, legend cell align=left, align=left, draw=none, fill=none, font=\footnotesize}
]
\addplot [color=mycolor1, line width=3.0pt, mark size=3.0pt, mark=o, mark options={solid, mycolor1}]
  table[row sep=crcr]{%
95	3.8804\\
767	14.6382\\
6143	16.1164\\
};
\addlegendentry{MLVSBM}

\addplot [color=mycolor2, line width=3.0pt, mark size=3.0pt, mark=+, mark options={solid, mycolor2}]
  table[row sep=crcr]{%
95	7.7508\\
767	11.11\\
6143	17.6321\\
};
\addlegendentry{MLVSMATCH3}

\addplot [color=mycolor3, line width=3.0pt, mark size=3.0pt, mark=asterisk, mark options={solid, mycolor3}]
  table[row sep=crcr]{%
95	6.0218\\
767	11.2393\\
6143	15.1177\\
};
\addlegendentry{MLVSMATCH4}

\addplot [color=mycolor4, line width=3.0pt, mark size=3.3pt, mark=triangle, mark options={solid, mycolor4}]
  table[row sep=crcr]{%
95	8.44837\\
767	14.6718\\
6143	31.8603\\
};
\addlegendentry{AlyaDefCG}

\end{axis}
\end{tikzpicture}}
\subfloat[\label{fig:weak2quater-bottom}]{\begin{tikzpicture}
\begin{axis}[%
width=0.25\columnwidth,
height=0.25\columnwidth,
at={(2.486in,1.101in)},
scale only axis,
xmin=191,
xmax=12287,
xtick={  191,  1535, 12287},
xticklabels={  191,  1535, 12287},
yticklabel style = {font=\footnotesize},
xticklabel style = {font=\footnotesize},
xmode=log,
xlabel style={font=\color{white!15!black}\footnotesize},
xlabel={number of MPI cores},
ymin=2.5079,
ymax=27.7146,
ylabel style={font=\color{white!15!black}\footnotesize},
axis background/.style={fill=white},
title style={font=\small},
legend style={at={(0.03,0.97)}, anchor=north west, legend cell align=left, align=left, draw=none, fill=none, font=\footnotesize}
]
\addplot [color=mycolor1, line width=3.0pt, mark size=3.0pt, mark=o, mark options={solid, mycolor1}]
  table[row sep=crcr]{%
191	2.5079\\
1535	7.1922\\
12287	12.375\\
};
\addlegendentry{MLVSBM}

\addplot [color=mycolor2, line width=3.0pt, mark size=3.0pt, mark=+, mark options={solid, mycolor2}]
  table[row sep=crcr]{%
191	4.3741\\
1535	8.3821\\
12287	17.8802\\
};
\addlegendentry{MLVSMATCH3}

\addplot [color=mycolor3, line width=3.0pt, mark size=3.0pt, mark=asterisk, mark options={solid, mycolor3}]
  table[row sep=crcr]{%
191	3.6659\\
1535	10.9611\\
12287	14.2742\\
};
\addlegendentry{MLVSMATCH4}

\addplot [color=mycolor4, line width=3.0pt, mark size=3.3pt, mark=triangle, mark options={solid, mycolor4}]
  table[row sep=crcr]{%
191	8.44837\\
1535	8.96248\\
12287	27.7146\\
};
\addlegendentry{AlyaDefCG}

\end{axis}
\end{tikzpicture}}%

%% file: weakspeedupofsolvemarenostrum.tikz
\definecolor{mycolor1}{rgb}{0.00000,0.44700,0.74100}%
\definecolor{mycolor2}{rgb}{0.85000,0.32500,0.09800}%
\definecolor{mycolor3}{rgb}{0.92900,0.69400,0.12500}%
\definecolor{mycolor4}{rgb}{0.49400,0.18400,0.55600}%
\subfloat[\label{weak2tris-top}]{\begin{tikzpicture}
\begin{axis}[%
width=0.25\columnwidth,
height=0.25\columnwidth,
at={(2.5in,7.1in)},
scale only axis,
xmin=47,
xmax=3071,
xtick={  47,  383, 3071},
xticklabels={  47,  383, 3071},
yticklabel style = {font=\footnotesize},
xticklabel style = {font=\footnotesize},
xmode=log,
xlabel style={font=\color{white!15!black}\footnotesize},
xlabel={number of MPI cores},
%ymode=log,
ymin=1,
ymax=45.1659403631988,
yminorticks=false,
ylabel style={font=\color{white!15!black}\footnotesize},
axis background/.style={fill=white},
title style={font=\footnotesize},
legend style={at={(0.97,0.33)}, anchor=south east, legend cell align=left, align=left, draw=none, fill=none, font=\footnotesize}
]
\addplot [color=mycolor1, line width=3.0pt, mark size=3.0pt, mark=o, mark options={solid, mycolor1}]
  table[row sep=crcr]{%
47	1\\
383	2.69265110808138\\
3071	25.6215030286329\\
};
\addlegendentry{MLVSBM}

\addplot [color=mycolor2, line width=3.0pt, mark size=3.0pt, mark=+, mark options={solid, mycolor2}]
  table[row sep=crcr]{%
47	1\\
383	5.36218243514756\\
3071	45.1659403631988\\
};
\addlegendentry{MLVSMATCH3}

\addplot [color=mycolor3, line width=3.0pt, mark size=3.0pt, mark=asterisk, mark options={solid, mycolor3}]
  table[row sep=crcr]{%
47	1\\
383	4.55109811915843\\
3071	36.8510141316619\\
};
\addlegendentry{MLVSMATCH4}

\addplot [color=mycolor4, line width=3.0pt, mark size=3.3pt, mark=triangle, mark options={solid, mycolor4}]
  table[row sep=crcr]{%
47	1\\
383	4.51526443223334\\
3071	23.2558909204123\\
};
\addlegendentry{AlyaDefCG}

\end{axis}
\end{tikzpicture}}
\subfloat[\label{weak2tris-middle}]{
\begin{tikzpicture}
\begin{axis}[%
width=0.25\columnwidth,
height=0.25\columnwidth,
at={(2.5in,7.1in)},
scale only axis,
xmin=95,
xmax=6143,
xtick={  95,  767, 6143},
xticklabels={  95,  767, 6143},
yticklabel style = {font=\footnotesize},
xticklabel style = {font=\footnotesize},
xmode=log,
xlabel style={font=\color{white!15!black}\footnotesize},
xlabel={number of MPI cores},
%ymode=log,
ymin=1,
ymax=28.1334157587582,
yminorticks=false,
ylabel style={font=\color{white!15!black}\footnotesize},
axis background/.style={fill=white},
title style={font=\footnotesize},
legend style={at={(0.97,0.33)}, anchor=south east, legend cell align=left, align=left, draw=none, fill=none, font=\footnotesize}
]
\addplot [color=mycolor1, line width=3.0pt, mark size=3.0pt, mark=o, mark options={solid, mycolor1}]
  table[row sep=crcr]{%
95	1\\
767	2.12069790001503\\
6143	15.4094959172024\\
};
\addlegendentry{MLVSBM}

\addplot [color=mycolor2, line width=3.0pt, mark size=3.0pt, mark=+, mark options={solid, mycolor2}]
  table[row sep=crcr]{%
95	1\\
767	5.58113411341134\\
6143	28.1334157587582\\
};
\addlegendentry{MLVSMATCH3}

\addplot [color=mycolor3, line width=3.0pt, mark size=3.0pt, mark=asterisk, mark options={solid, mycolor3}]
  table[row sep=crcr]{%
95	1\\
767	4.28624558468944\\
6143	25.4929784292584\\
};
\addlegendentry{MLVSMATCH4}

\addplot [color=mycolor4, line width=3.0pt, mark size=3.3pt, mark=triangle, mark options={solid, mycolor4}]
  table[row sep=crcr]{%
95	1\\
767	4.60658951185267\\
6143	16.9708282721757\\
};
\addlegendentry{AlyaDefCG}

\end{axis}
\end{tikzpicture}}
\subfloat[\label{weak2tris-bottom}]{\begin{tikzpicture}
\begin{axis}[%
width=0.25\columnwidth,
height=0.25\columnwidth,
at={(2.5in,7.1in)},
scale only axis,
xmin=191,
xmax=12287,
xtick={  191,  1535, 12287},
xticklabels={  191,  1535, 12287},
yticklabel style = {font=\footnotesize},
xticklabel style = {font=\footnotesize},
xmode=log,
xlabel style={font=\color{white!15!black}\footnotesize},
xlabel={number of MPI cores},
%ymode=log,
ymin=1,
ymax=19.5094166973364,
yminorticks=false,
ylabel style={font=\color{white!15!black}\footnotesize},
axis background/.style={fill=white},
title style={font=\footnotesize},
legend style={at={(0.97,0.33)}, anchor=south east, legend cell align=left, align=left, draw=none, fill=none, font=\footnotesize}
]
\addplot [color=mycolor1, line width=3.0pt, mark size=3.0pt, mark=o, mark options={solid, mycolor1}]
  table[row sep=crcr]{%
191	1\\
1535	2.78957759795334\\
12287	12.9701494949495\\
};
\addlegendentry{MLVSBM}

\addplot [color=mycolor2, line width=3.0pt, mark size=3.0pt, mark=+, mark options={solid, mycolor2}]
  table[row sep=crcr]{%
191	1\\
1535	4.17470562269598\\
12287	15.6565586514692\\
};
\addlegendentry{MLVSMATCH3}

\addplot [color=mycolor3, line width=3.0pt, mark size=3.0pt, mark=asterisk, mark options={solid, mycolor3}]
  table[row sep=crcr]{%
191	1\\
1535	2.67557088248442\\
12287	16.436479802721\\
};
\addlegendentry{MLVSMATCH4}

\addplot [color=mycolor4, line width=3.0pt, mark size=3.3pt, mark=triangle, mark options={solid, mycolor4}]
  table[row sep=crcr]{%
191	1\\
1535	7.54110023118601\\
12287	19.5094166973364\\
};
\addlegendentry{AlyaDefCG}

\end{axis}
\end{tikzpicture}}

%% file: weakavgiterationjuwels.tikz
\definecolor{mycolor1}{rgb}{0.00000,0.44700,0.74100}%
\definecolor{mycolor2}{rgb}{0.85000,0.32500,0.09800}%
\definecolor{mycolor3}{rgb}{0.92900,0.69400,0.12500}%
\definecolor{mycolor4}{rgb}{0.49400,0.18400,0.55600}%
\begin{tikzpicture}

\begin{axis}[%
width=0.5\columnwidth,
height=1.6in,
at={(0.758in,0.481in)},
scale only axis,
xmin=0,
xmax=3,
xtick={0,1,2,3,4},
xticklabels={{45},{367},{2943},{23551},{}},
xlabel style={font=\color{white!15!black}},
xlabel={Number of MPI cores},
ymin=1,
ymax=350,
ytick={4,100,300},
yticklabels={4,100,300},
ymode=log,
ylabel style={font=\color{white!15!black}},
ylabel={Average Iterations},
axis background/.style={fill=white},
legend style={legend cell align=left, align=left, draw=none, fill=none}
]
\addplot [color=mycolor4, line width=3.0pt, mark size=3.3pt, mark=triangle, mark options={solid, mycolor4}]
  table[row sep=crcr]{%
0	132.619047619048\\
1	95.2380952380952\\
2	149.142857142857\\
3	330.952380952381\\
};
\addlegendentry{AlyaCG}

\addplot [color=mycolor3, line width=3.0pt, mark size=3.0pt, mark=asterisk, mark options={solid, mycolor3}]
  table[row sep=crcr]{%
0	6\\
1	6\\
2	4.80952380952381\\
3	4\\
};
\addlegendentry{MLVSMATCH4}

\end{axis}
\end{tikzpicture}%

%% file: weaktotalsolvetimejuwels.tikz
\definecolor{mycolor1}{rgb}{0.00000,0.44700,0.74100}%
\definecolor{mycolor2}{rgb}{0.85000,0.32500,0.09800}%
\definecolor{mycolor3}{rgb}{0.92900,0.69400,0.12500}%
\definecolor{mycolor4}{rgb}{0.49400,0.18400,0.55600}%
\begin{tikzpicture}

\begin{axis}[%
width=0.35\columnwidth,
height=1.5in,
at={(0.758in,0.481in)},
scale only axis,
xmin=0,
xmax=3,
xtick={0,1,2,3,4},
xticklabels={{45},{367},{2943},{23551},{}},
xlabel style={font=\color{white!15!black}},
xlabel={Number of MPI cores},
ymin=10,
ymax=100,
axis background/.style={fill=white},
legend style={legend cell align=left, align=left, draw=none, fill=none}
]
\addplot [color=mycolor4, line width=3.0pt, mark size=3.3pt, mark=triangle, mark options={solid, mycolor4}]
  table[row sep=crcr]{%
0	32.83255943\\
1	25.15684273\\
2	40.63196319\\
3	95.05365349\\
};
\addlegendentry{AlyaCG}

\addplot [color=mycolor3, line width=3.0pt, mark size=3.0pt, mark=asterisk, mark options={solid, mycolor3}]
  table[row sep=crcr]{%
0	16.79627197\\
1	16.847301\\
2	15.06735504\\
3	13.42745439\\
};
\addlegendentry{MLVSMATCH4}

\end{axis}
\end{tikzpicture}%

%% file: weakspeedupofsolvejuwels.tikz
\definecolor{mycolor1}{rgb}{0.00000,0.44700,0.74100}%
\definecolor{mycolor2}{rgb}{0.85000,0.32500,0.09800}%
\definecolor{mycolor3}{rgb}{0.92900,0.69400,0.12500}%
\definecolor{mycolor4}{rgb}{0.49400,0.18400,0.55600}%
\begin{tikzpicture}

\begin{axis}[%
width=0.35\columnwidth,
height=1.5in,
at={(0.758in,0.481in)},
scale only axis,
xmin=0,
xmax=3,
xtick={0,1,2,3,4},
xticklabels={{45},{367},{2943},{23551},{}},
xlabel style={font=\color{white!15!black}},
xlabel={Number of MPI cores},
ymin=1,
ymax=700,
ytick={1,10,100,700},
yticklabels={1,10,100,700},
ymode=log,
axis background/.style={fill=white},
legend style={legend cell align=left, align=left, draw=none, fill=none}
]
\addplot [color=mycolor4, line width=3.0pt, mark size=3.3pt, mark=triangle, mark options={solid, mycolor4}]
  table[row sep=crcr]{%
0	1\\
1	10.4409157484128\\
2	51.7150449682714\\
3	176.850334636832\\
};
\addlegendentry{AlyaCG}

\addplot [color=mycolor3, line width=3.0pt, mark size=3.0pt, mark=asterisk, mark options={solid, mycolor3}]
  table[row sep=crcr]{%
0	1\\
1	7.97576868603464\\
2	71.3437363907767\\
3	640.455815291732\\
};
\addlegendentry{MLVSMATCH4}

\end{axis}
\end{tikzpicture}%